\documentclass[reqno,12pt]{amsart}  

\usepackage{gensymb}
 
 \usepackage{amsmath}
\usepackage{amssymb} 
\usepackage{mathtools}   
\usepackage{graphicx}
\usepackage{graphics}  
\usepackage{color}
\usepackage{verbatim} 
\usepackage[normalem]{ulem} 
\usepackage[mathscr]{eucal}
\usepackage{upgreek}
\usepackage{enumerate}

\usepackage[titletoc]{appendix}
\numberwithin{equation}{section}  
 
\usepackage[refpage,noprefix]{nomencl}
\usepackage{nomencl} 

\makenomenclature

\newtheorem{theorem}{Theorem}[section] 
\newtheorem{lemma}[theorem]{Lemma} 
\newtheorem{proposition}[theorem] {Proposition} 
 
\newtheorem{remark}[theorem]  {Remark} 
\newtheorem{definition}[theorem] {Definition}

\theoremstyle{definition}

 \usepackage{float}
\floatstyle{boxed}
\restylefloat{figure}

\DeclareMathAlphabet{\mathpzc}{OT1}{pzc}{m}{it}

\DeclarePairedDelimiter{\norm}{\lVert}{\rVert}

\renewcommand{\L} {\Lambda} %

\def\d{\delta} 
 
\newcommand{\eps}{\varepsilon}

\def\t{\tau}

\newfam\Bbbfam 
\font\tenBbb=msbm10 
\font\sevenBbb=msbm7 
\font\fiveBbb=msbm5 
\textfont\Bbbfam=\tenBbb 
\scriptfont\Bbbfam=\sevenBbb 
\scriptscriptfont\Bbbfam=\fiveBbb

\newcommand{\R}     {\mathbb{R}} 
\newcommand{\Z}     {\mathbb{Z}} 
\newcommand{\N}     {\mathbb{N}} 
\renewcommand{\P}   {\mathbb{P}} 
 
\newcommand{\E}     {\mathbb{E}}

 \newcommand{\floor}[1]{\left\lfloor #1 \right\rfloor}

\def\1{{\mathchoice {1\mskip-4mu\mathrm l}      
{1\mskip-4mu\mathrm l} 
{1\mskip-4.5mu\mathrm l} {1\mskip-5mu\mathrm l}}} 
\newcommand{\ssup}[1] {{\scriptscriptstyle{({#1}})}} 
\def\comment#1{} 
\newtheoremstyle{thm}{2ex}{2ex}{\itshape\rmfamily}{} 
{\bfseries\rmfamily}{}{1.7ex}{} 
 
\newtheoremstyle{rem}{1.3ex}{1.3ex}{\rmfamily}{} 
{\itshape\rmfamily}{}{1.5ex}{} 
 
\newenvironment{proofsect}[1] 
{\vskip0.1cm\noindent{\scshape #1.}\hskip0.5cm}

\newcommand{\ba}{\boldsymbol{a}}
\newcommand{\bb} {\boldsymbol{b}}

\newcommand{\br}{\boldsymbol{r}}

\newcommand{\zero} {\boldsymbol{0}}

\newcommand{\Ical}   {{\mathcal I }}

\newcommand{\Ascr} {\mathscr{A}}

\newcommand{\Cscr} {\mathscr{C}}
\newcommand{\Dscr}{\mathscr{D}}
\newcommand{\Escr}{\mathscr{E}}

\newcommand{\Hscr} {\mathscr{H}}

\newcommand{\Kscr}{\mathscr{K}}

\newcommand{\Mscr}{\mathscr{M}}
\newcommand{\Nscr}{\mathscr{N}}
\newcommand{\Oscr}{\mathscr{O}}
\newcommand{\Pscr}{\mathscr{P}}

\newcommand{\Sscr}{\mathscr{S}}

\newcommand{\gNr}{\gamma_{N,\eps}^{\br}}
\newcommand{\ZNe}{Z_{N,\eps}}
\newcommand{\gNll}{\gamma_{[-1,p_*+1]}^{\br}}
\newcommand{\ZNl}{Z_{[-1, p_*+1]}}

\newcommand{\ZNr}{Z_{[p^*, N+1]}}
\newcommand{\ZNrr}{Z_{[p^*-1, N+1]}}
\newcommand{\gNP}{\gamma_{[p_*, p^*] \setminus \Pscr}^{\zero}}
\newcommand{\ZNP}{Z_{[p_*, p^*] \setminus \Pscr} }
\newcommand{\gNPP}{\gamma_{[-1, N+1] \setminus \Pscr}^{\br}}
\newcommand{\ZNPP}{Z_{[-1, N+1] \setminus \Pscr}}

 \newcommand{\ex}{{\rm e}} 
 
\renewcommand{\d}{{\rm d}}

\newcommand{\sign}{{\operatorname {sign}\,}}
\newcommand{\dist}{{\operatorname {dist}}}

\newcommand{\Exp}{\mathscr{E}\kern-0.2mm{\operatorname{xp}}}
\newcommand{\Log}{\mathscr{L}\kern-0.2mm{\operatorname{og}}}

\renewcommand{\emptyset} {\varnothing}

\makeindex

\begin{document}

\title[\hfill Large deviations\hfill]
{Sample path large deviations for Laplacian models in $(1+1)$-dimensions}


\author{Stefan Adams, Alexander Kister, and  Hendrik Weber}
\address{Mathematics Institute, University of Warwick, Coventry CV4 7AL, United Kingdom}
\email{S.Adams@warwick.ac.uk}

\thanks{}
  

 


 \maketitle 
\begin{center}

\end{center}

 \begin{abstract} 
 We study scaling limits of a Laplacian pinning model in $(1+1)$  dimension and derive sample path large deviations for the profile height function. The model is given by a Gaussian integrated random walk (or a Gaussian  integrated random walk bridge) perturbed by an attractive force towards the zero-level.  We study in detail the behaviour of the rate function and show that   it can  admit  up to  five minimisers depending on the choices of pinning strength and boundary conditions.  This study complements corresponding large deviation results for Gaussian gradient systems with pinning  in $ (1+1) $-dimension (\cite{FS04})  in $(1+d) $-dimension (\cite{BFO}), and recently in higher dimensions in \cite{BCF}. 
\end{abstract}

\section{Introduction and large deviation results} 
\subsection{The models}
We are going to study models for  $(1+1)$-dimensional random fields. These models are defined in terms of the potential, a measurable function $ V\colon\R\to\R\cup\{+\infty\} $ such that $ x\mapsto \exp(-V(x)) $ is bounded and continuous and that 
$$ 
\int_\R\ex^{-V(x)}\,\d x<\infty \,\mbox{  and  } \, \int_\R x^2\ex^{-V(x)}\,\d x=:\sigma^2<\infty\,\mbox{ and } \int_\R x\ex^{-V(x)}\,\d x=0.
$$
For most of the article we consider the Gaussian case $V(x) = \frac12 x^2$.
Given the potential $V$, we define a Hamiltonian $ \Hscr_{[\ell,r]} (\phi)$, defined for $ \ell,r \in\Z$, with $ r-\ell\ge 2$, and for $ \phi\colon\{\ell,\ell+1,\ldots,r-1,r\} \to\R $ by
\begin{equation}\label{Hamiltonian}
\Hscr_{[\ell,r]}(\phi):=\sum_{k=\ell+1}^{r-1}V(\Delta\phi_k),
\end{equation}
where $ \Delta $ denotes the discrete Laplacian, $\Delta\phi_k=\phi_{k+1}+\phi_{k-1}-2\phi_k $. Our  pinning models are then given by the probability measures 
\begin{equation}\label{gibbsdist}
\begin{aligned}
\gamma_{N,\eps}^{\psi}(\d\phi)&=\frac{1}{\ZNe(\psi)}\ex^{-\Hscr_{[-1,N+1]}(\phi)}\prod_{k=1}^{N-1}(\eps\delta_0(\d\phi_k)+\d\phi_k)\prod_{k\in\{-1,0,N,N+1\}}\delta_{\psi_k}(\d\phi_k),\\
\gamma_{N,\eps}^{\psi_f}(\d\phi)&=\frac{1}{\ZNe(\psi_f)}\ex^{-\Hscr_{[-1,N+1]}(\phi)}\prod_{k=1}^{N+1}(\eps\delta_0(\d\phi_k)+\d\phi_k)\prod_{k\in\{-1,0\}}\delta_{\psi_k}(\d\phi_k),
\end{aligned}
\end{equation}
where $ N\ge 2 $ is an integer, $ \eps\ge 0 $ is the pinning strength, $\d\phi_k $  is the Lebesgue measure on $ \R$, $\delta_0 $ is the Dirac mass at zero, where $\psi\in\R^\Z $ is a given boundary condition and $ \ZNe(\psi) $ (resp. $ \ZNe(\psi_f) $) is the normalisation, which is usually called partition function. 

The measures given  in \eqref{gibbsdist}  are $ (1+1) $-dimensional models for a linear chain of length $N$ which is attracted to the defect line, the $x$-axis. The parameter $ \eps \ge 0 $  tunes the strength of the attraction and one wishes to understand its effect on the field, in the large $N$ limit. The models with $ \eps=0 $ have no pinning reward at all and are thus free Laplacian models. By``$(1+1)$-dimensional" we mean that the configurations of the chain are given by graphs $ \{(k,\phi_k)\}_{-1\le k\le N+1} $.   Models with Laplacian interaction have been studied in the Physics literature in the context of semiflexible polymers, c.f. \cite{BLL,HV09}, or in the context of deforming rods in space, cf. \cite{Antman}.

The basic properties of the models  were investigated in the two papers \cite{CD08,CD09}, to which we refer for a detailed discussion and  for a survey of the literature. In particular, it  was shown in \cite{CD08} that there is  a critical value $ \eps_{\rm c}\in (0,\infty) $  that determines a  phase transitions between a \textit{delocalised regime} ($\eps<\eps_{\rm c}$), in which the reward is essentially ineffective, and a \textit{localised regime} ($\eps>\eps_{\rm c}$), in which  the reward has a macroscopic effect on the field. For more details see Section~\ref{Sec-LDPpinning} below. 
In the present paper we derive large deviation principles for the macroscopic empirical profile distributed under the measures in \eqref{gibbsdist} (Section~\ref{Sec-LDP}). The corresponding large deviation results  for the gradient models, where the Hamiltonian is a function of the discrete gradient of the field instead of the discrete Laplacian, have been derived  in \cite{FS04} for Gaussian random walk bridges in $ \R $ and for Gaussian random walks and bridges in higher dimensions in \cite{BFO}. In \cite{FO10} large deviations for general non-Gaussian random walks in $ \R^d,d\ge 1 $, have were analysed, and in \cite{BCF} gradient model in higher (lattice) dimensions were introduced. 

A common feature of all these gradient models is, that typical fluctuations are observed on scale $ \sqrt{N} $ and that large deviation results can be obtained on a  linear scale in $N$; for more details see \cite{Funaki}. In contrast  in the Laplacian case  the scale for the scaling limits is $ N^{3/2} $ as already observed in Sinai's work \cite{Sinai} on integrated random walks and proved in the specific context of our models by Caravenna and Deuschel in \cite{CD09}. In this article we derive large deviations principles on scale  $N^2$.
Beyond the different scaling, a major technical difference between the Laplacian case and the gradient case is the fact that the Markov property which features prominently in the large deviation proofs in the gradient case \cite{FS04,BFO} is not directly available in the Laplacian case. To overcome this difficulty  we introduce a correction technique and replace ``single zeros'' of the profile by ``double zeros'' which then allows us to write the distribution over disjoint intervals separated by a double zero as  the product of independent distributions over the disjoint intervals. 

%
 
Our second major result (given in Section~\ref{Sec-Min}) is a complete analysis of the rate functions. We are particularly interested  in the critical situation where  more than one minimiser of the  rate function exists. 
The variational problem for our rate functions shows a much richer structure of minimisers than the rate function for corresponding gradient models in \cite{BFO}. Without pinning there is a unique bi-harmonic function minimising the macroscopic bi-Laplacian energy (see Appendix~\ref{AppA}). Once the pinning reward is switched on, the integrated random walk (scaled random field) has essentially two different strategies to pick up reward. One strategy is to start picking up the reward earlier despite the energy involved to bend to the zero line with speed zero and the other strategy is to cross the zero level producing a longer bend before turning to the zero level and picking up reward. The choice of pinning strength and boundary conditions determines which of these strategies is favoured by the rate function.

In Section~\ref{Sec-LDP} we present the large deviation results which are proved in Section~\ref{Sec-PLDP}. The results of the variational analysis  are given in Section~\ref{Sec-Min} and  their proofs are given in Section~\ref{Sec-minproof}. We include an Appendix~\ref{AppA} where some basic facts about the bi-harmonic equation and bi-harmonic functions along with convergence statements for the discrete bi-Laplacian are provided. In Appendix~\ref{AppB} we collect  some well-known facts about  partition functions of Gaussian integrated random walks.

\subsection{Sample path large deviations}\label{Sec-LDP}

 \subsubsection{Empirical profile} Let $ h_N=\{h_N(t)\colon t\in[0,1]\} $ be the macroscopic empirical profile determined from the microscopic height function $ \phi $ under the proper scaling. More precisely, define $h_N$  as a linear interpolation of $ (h_N(k/N)=\phi_k/N^2)_{k\in\L_N}$ with $ \L_N=\{-1,0,\ldots, N,N+1\} $ by
  \begin{equation}\label{def:hN}
 h_N(t)=\frac{\floor{Nt}-Nt+1}{N^2}\phi_{\floor{Nt}}+\frac{Nt-\floor{Nt}}{N^2}\phi_{\floor{Nt}+1},\quad t\in[0,1].
 \end{equation}
 We study $ h_N $ distributed under the  measures given in \eqref{gibbsdist} endowed with a suitable boundary conditions $ \psi^{\ssup{N}} $. In the case of Dirichlet boundary conditions, we fix parameters  $ a,\alpha,b,\beta\in\R $ and then define the microscopic boundary conditions as
\begin{equation}\label{bc}
 \psi^{\ssup{N}}(x)=\begin{cases} aN^2-\alpha N &\mbox{ if } x=-1,\\
 aN^2 & \mbox{ if } x=0,\\
 bN^2 & \mbox{ if } x=N,\\
 bN^2+\beta N & \mbox{ if } x=N+1,\\
 0 & \mbox{ otherwise}.\end{cases}
 \end{equation}
 On the macroscopic scale this choice corresponds to fixing 
%
 $ h_N(0)=a $ and $  h_N(1)=b $ as well the discrete derivatives 
 \begin{equation*}
  \dot{h}_N(0)=\frac{\psi^{\ssup{N}}(0)-\psi^{\ssup{N}}(-1)}{N}=\alpha \quad \text{and} \quad  \dot{h}_N(1)=\frac{\psi^{\ssup{N}}(N+1)-\psi^{\ssup{N}}(N)}{N}=\beta .
  \end{equation*} 
In the case of free boundary conditions on the right side
we only specify the boundary in $ x=-1 $ and $x=0 $, and write  $ \psi_f^{\ssup{N}}(-1)=aN^2-\alpha N $ and $ \psi_f^{\ssup{N}}(0)=aN^2 $, see \eqref{gibbsdist}.  
We write $ \br=(a,\alpha,b,\beta) $ to specify our choice of boundary conditions $ \psi^{\ssup{N}}$ in the Dirichlet case and $ \ba=(a,\alpha) $ for the mixed Dirichlet and free boundary case.

We denote the Gibbs distributions with $ \eps=0 $ (no pinning) by $ \gamma_N^{\br} $ for Dirichlet boundary conditions and by  $ \gamma_{N}^{\ba} $ for Dirichlet boundary conditions on the left  and free boundary conditions on the right  and their partition functions by $ Z_N(\br) $ and $ Z_N(\ba) $, respectively.  In Section~\ref{Sec-LDPnp} we study the large deviation principles  without pinning ($\eps=0$) for general integrated random walks with free boundary conditions on the right hand side and show that these results apply to  Gaussian integrated random walk bridges  as well. 
Our main large deviation result for the measures with pinning are then presented in Section~\ref{Sec-LDPpinning}. 
To state these results we introduce the spaces 
$$ H^2_{\br}=\{h\in H^2([0,1])\colon h(0)=a,h(1)=b,\dot{h}(0)=\alpha,\dot{h}(1)=\beta\} ,
$$ 
used for the case of Dirichlet boundary conditions and the space
$$ H^2_{\ba}=\{h\in H^2([0,1])\colon h(0)=a,\dot{h}(0)=\alpha\}
$$  
used in the case of free boundary conditions on the right. 
Here, $ H^2([0,1]) $ is the usual Sobolev space. We write $ \Cscr([0,1];\R) $ for the space of continuous functions on $ [0,1] $ equipped with the supremum norm.

\subsubsection{Large deviations for integrated random walks and Gaussian integrated random walk bridges}\label{Sec-LDPnp}

We recall the integrated random walk representation in Proposition~2.2 of \cite{CD08}. Let $ (X_k)_{k\in\N} $ be a sequence of independent and identically distributed random variables, with marginal laws $ X_1\sim\exp(-V(x))\d x $, and $ (Y_n)_{n\in\N_0} $  the corresponding random walk with initial condition $ Y_0=\alpha N $ and $ Y_n=\alpha N +X_1+\cdots + X_n $. The integrated random walk is denoted by $ (Z_n)_{n\in\N_0} $ with $ Z_0=aN^2 $ and $ Z_n=aN^2 +Y_1+\cdots + Y_n $. We denote $ \P^{\ba} $ the probability distribution of the above defined processes.
Then the following holds.
\begin{proposition}[\textbf{\cite{CD08}}]\label{P:CD}
The pinning free model   $ \gamma_{N}^{\br} $ ($\eps = 0$) is the law of the vector $ (Z_1,\ldots, Z_{N-1}) $ under the measure $ \P^{\br}(\cdot):=\P^{\ba}(\cdot |Z_N=bN^2,Z_{N+1}=bN^2+\beta N) $. The partition function $ Z_N(\br) $ is the value at $ (\beta N, bN^2+\beta N) $ of the density of the vector $ (Y_{N+1},Z_{N+1}) $ under the law $ \P^{\ba} $. The model $ \gamma_N^{\ba} $ coincides with the integrated random walk $ \P^{\ba} $.

\end{proposition}

The first part of the following result is the generalisation of Mogulskii's theorem \cite{Mog} from random walks to integrated random walks whereas its second part is the generalisation to Gaussian integrated random walk bridges.

\begin{theorem}\label{LDP-nopinning}
\begin{enumerate}

\item[(a)] Let $ V$ be any potential of the form above such that $ \L(\lambda):=\log\E[\ex^{\langle \lambda,X_1\rangle}]<\infty $ for all $ \lambda\in\R $, then the following holds. The large deviation principle (LDP) holds for $ h_N $ under $ \gamma_{N}^{\ba} $  on the space $ \Cscr([0,1];\R) $ as $ N\to\infty $ with speed $N$ and the unnormalised good rate function  $ \Escr_{f} $ of the form:
\begin{equation}
\Escr_f(h)=\begin{cases} \int_0^1\,\L^*(\ddot{h}(t))\,\d t, & \mbox{ if } h\in H^2_{\ba},\\ +\infty & \mbox{ otherwise}.\end{cases}
\end{equation}   Here $ \L^* $ denotes the Fenchel-Legendre transform of $ \L $.

\item[(b)] For $ V(\eta)=\frac{1}{2}\eta^2 $ the following holds.
The large deviation principle (LDP) holds for $ h_N $ under $ \gamma_N^{\br}  $ on the spaces $ \Cscr([0,1];\R) $  as $ N\to\infty $ with speed $N$ and the  unormalised good rate function  $ \Escr $ of the form:
\begin{equation}\label{e:onesix}
\Escr(h)=\begin{cases} \frac{1}{2}\int_0^1\ddot{h}^2(t)\,\d t    , &\mbox{ if } h\in H^2_{\br},\\
+\infty & \mbox{ otherwise}.
\end{cases}
\end{equation}

\end{enumerate}
\end{theorem}

\begin{remark}
\begin{enumerate}
\item[(a)] The rate functions in both cases are obtained from the unnormalised rate functions by $ I_f^0(h)=\Escr_f(h)-\inf_{g\in H^2_{\ba}}\Escr_f(g) $ for  general integrated random walks with potential $ V$ respectively by
$ I^0(h)=\Escr(h)-\inf_{g\in H^2_{\br}}\Escr(g) $ for  Gaussian integrated random walk bridges.

\item[(b)] We believe that the large deviation in Theorem~\ref{LDP-nopinning}(b)  holds for general potentials $ V$ as in (a)  as well. For the Gaussian integrated random walk bridges there exist   explicit formulae for  the distribution, see \cite{GSV}.  Our main result concerns the large deviations for the pinning model for Gaussian integrated random walk bridges. General integrated random walk bridges will require  different techniques. 

\end{enumerate}
\end{remark}
\medskip

\medskip

\subsubsection{Large deviations for pinning models}\label{Sec-LDPpinning}

The large deviation principle for the pinning models gets an additional term for  the rate function. Recall that the logarithm of the partition function is the free energy. Difference of  the free energies with pinning and without pinning for zero boundary conditions $ (\br=\zero$) will be an important ingredient in our rate functions. We define  $ \tau(\eps) $ as the thermodynamic limit of the logarithm of the quotient of the partition function with pinning and the partition function without pinning (both with zero boundary condition),
\begin{equation}\label{freeenergies}
\tau(\eps)=\lim_{N\to\infty}\frac{1}{N}\log\frac{\ZNe(\zero)}{Z_N(\zero)}.
\end{equation}
 The existence of the limit in \eqref{freeenergies} and its  properties have been derived by Caravenna and Deuschel in \cite{CD08}, we summarise their result in the following proposition.

\begin{proposition}[\cite{CD08}]\label{existenceOfLimit}
The limit in \eqref{freeenergies} exist for every $ \eps\ge 0$. Furthermore, there exists  $ \eps_{\rm c} \in (0,\infty) $ such that $ \tau(\eps)=0 $ for $ \eps\in[0,\eps_{\rm c}] $, while $ 0<\tau(\eps)<\infty $ for $ \eps\in (\eps_{\rm c},\infty) $, and as $ \eps\to\infty $,
$$
\tau(\eps)=\log\eps(1-o(1)).
$$ Moreover the function $ \tau $ is real analytic on $ (\eps_{\rm c},\infty) $.
\end{proposition} 
 We have the following sample path large deviation principles for $ h_N $ under $ \gamma_{N,\eps}^{\br} $ and  $ \gamma_{N,\eps}^{\ba}$, respectively.  The unnormalised rate functions denoted by $ \Sigma^\eps $ and  $ \Sigma^\eps_f $ are of the form
 \begin{equation}\label{ratefctSigma}
 \Sigma^\eps(h)=\frac{1}{2}\int_0^1 \ddot{h}^2(t)\,\d t-\tau(\eps)|\{t\in[0,1]\colon h(t)=0\}|,
 \end{equation}
 for $h\in H=H_{\br}^2 $ and $ H=H_{\ba}^{2} $, respectively. Here $ |\cdot| $ stands for the Lebesgue measure. 
  \begin{theorem}\label{THM-main} Let $ V(\eta)=\frac{1}{2}\eta^2 $.  The LDP holds for $h_N $ under $ \gamma_N=\gamma_{N,\eps}^{\br},\gamma_{N,\eps}^{\ba} $ respectively on the space $ \Cscr([0,1];\R) $ as $ N\to\infty $ with the speed $N$ and the good rate functions $ I=I^\eps $ and $  I=I^\eps_f  $ of the form:
 \begin{equation}
 I(h)=\begin{cases} \Sigma(h)-\inf_{h\in H}\{\Sigma(h)\} ,&\mbox{ if } h\in H,\\+\infty & \mbox{ otherwise},\end{cases}
 \end{equation}
 with $ \Sigma=  \Sigma^\eps $ and $ \Sigma=\Sigma^\eps_f $ respectively, and $ H=H_{\br}^2 $ respectively $H= H_{\ba}^{2}$. Namely, for every open set $ \Oscr $ and every closed set $ \Kscr $ of $ \Cscr([0,1];\R)$  equipped with the uniform topology, we have that
 \begin{align}
 \liminf_{N\to\infty}\frac{1}{N}\log\gamma_N(h_N\in \Oscr)&\ge -\inf_{h\in\Oscr} I(h), \notag\\
 \limsup_{N\to\infty}\frac{1}{N}\log\gamma_N(h_N\in\Kscr)&\le -\inf_{h\in\Kscr} I(h),\label{LDPBounds}
 \end{align}
 in each of two situations.
 \end{theorem}

 \bigskip

As the limit $ \tau(\eps) $ of the difference of the free energies appears in our rate functions it is worth pointing out that this has a direct translation in terms of  path properties of the field, see \cite{CD08}. This is the microscopic counterpart of the effect of the reward term in our pinning rate functions. Defining the contact number $ \ell_N $ by
$$
\ell_N:=\#\{k\in\{1,\ldots,N\}\colon \phi_k=0\},
$$
we can easily obtain that for $ \eps>0 $ (see \cite{CD08}),
$$
D_N(\eps):=\E_{\gamma_{N,\eps}^{\zero}}\big[\ell_N/N\big]\to\eps\tau^\prime(\eps) \; \mbox{ as } N\to\infty.
$$
This gives the following paths properties. When $ \eps>\eps_{\rm c} $, then $ D_N(\eps)\to D(\eps)>0 $ as $ N\to\infty $, and the mean contact density is non-vanishing leading to localisation of the field (integrated random walk respectively integrated random walk bridge). For the other case, $ \eps<\eps_{\rm c} $, we get $ D_N(\eps)\to 0 $ as $ N\to\infty $ and thus the contact density is vanishing in the  thermodynamic limits leading to de-localisation.
\section{Minimisers of the rate functions}\label{Sec-Min}

We are concerned with  the set $ \Mscr^\eps $ of the minimiser of the unnormalised rate functions in \eqref{ratefctSigma} for our pinning LDPs. Any minimiser of \eqref{ratefctSigma} is a zero of the corresponding rate function in Theorem~\ref{THM-main}.  We let $ h_{\br}^* \in H_{\br}^2 $  be the unique minimiser of the energy $ \Escr $ defined in \eqref{e:onesix} (see Proposition~\ref{uniquemin}), that is, $ \Escr(h)=1/2\int_0^1\ddot{h}^2(t)\,\d t $ is the energy of the bi-Laplacian in dimension one. 
  For any interval $ I\subset[0,1] $ we  let $ h_{\br}^{*,I} \in H_{\br}^2(I)$, where the boundary conditions apply to the boundaries of  $I$,  be the unique minimiser of $ \Escr^I(h)=\frac{1}{2}\int_I \,\ddot{h}^2(t)\,\d t $, and we sometimes write $ \ba,\bb $ for the boundary condition $ \br $ with $ \ba=(a,\alpha) $ and $ \bb=(b,\beta) $. Of major interest are the zero sets 
$$ \Nscr_h=\{t\in[0,1]\colon h(t)=0\}  \quad\mbox{ of any minimiser } h.
$$ 
In Section~\ref{Sec-Minfree} we study the minimiser for the case of Dirichlet boundary conditions on the left hand side and free boundary conditions on the right hand side, in Section~\ref{Sec-MinD} we summarise our findings for the Dirichlet boundary case on both the right hand and left hand side. In Section~\ref{Sec-minproof} we give the proofs for our statements.

\subsection{Free boundary conditions on the right hand side}\label{Sec-Minfree} 
We consider Dirichlet boundary conditions on the left hand side and the free boundary condition on the right side  only. 

Let $ \Mscr_f ^\eps$ denote  the set of minimiser of  $ \Sigma_f^\eps $.
\begin{proposition}\label{P:minfree}
For any boundary condition $ \ba=(a,\alpha) $ on the left hand side the set $ \Mscr_f^\eps $   of minimiser of $ \Sigma_f^\eps $  is a subset of
\begin{equation}\label{superset}
\{\overline{h}\}\cup\{h_\ell\colon \ell\in(0,1)\},
\end{equation}
where for any $ \ell\in(0,1) $ the functions $ h_\ell\in H_{\ba,f}^2 $   are given by
\begin{equation}
h_\ell(t)=\begin{cases} h_{(\ba,\zero)}^{*,(0,\ell)}(t) &,\mbox{ for } t\in [0,\ell),\\
0 &, \mbox{ for } t\in [\ell,1];\end{cases}
\end{equation}
and the function $ \overline{h}\in H_{\ba}^2 $  is the linear function $ \overline{h}(t)=a+\alpha t, t\in[0,1] $. 
\end{proposition}

\begin{figure}[h!]
\caption{$h_{\ell_1} $ for $a=1 $ and $ \alpha=-12,\tau=288,\ell_1=1/2(\sqrt{2}-1) $} \label{fig2}
\includegraphics[scale=0.42]{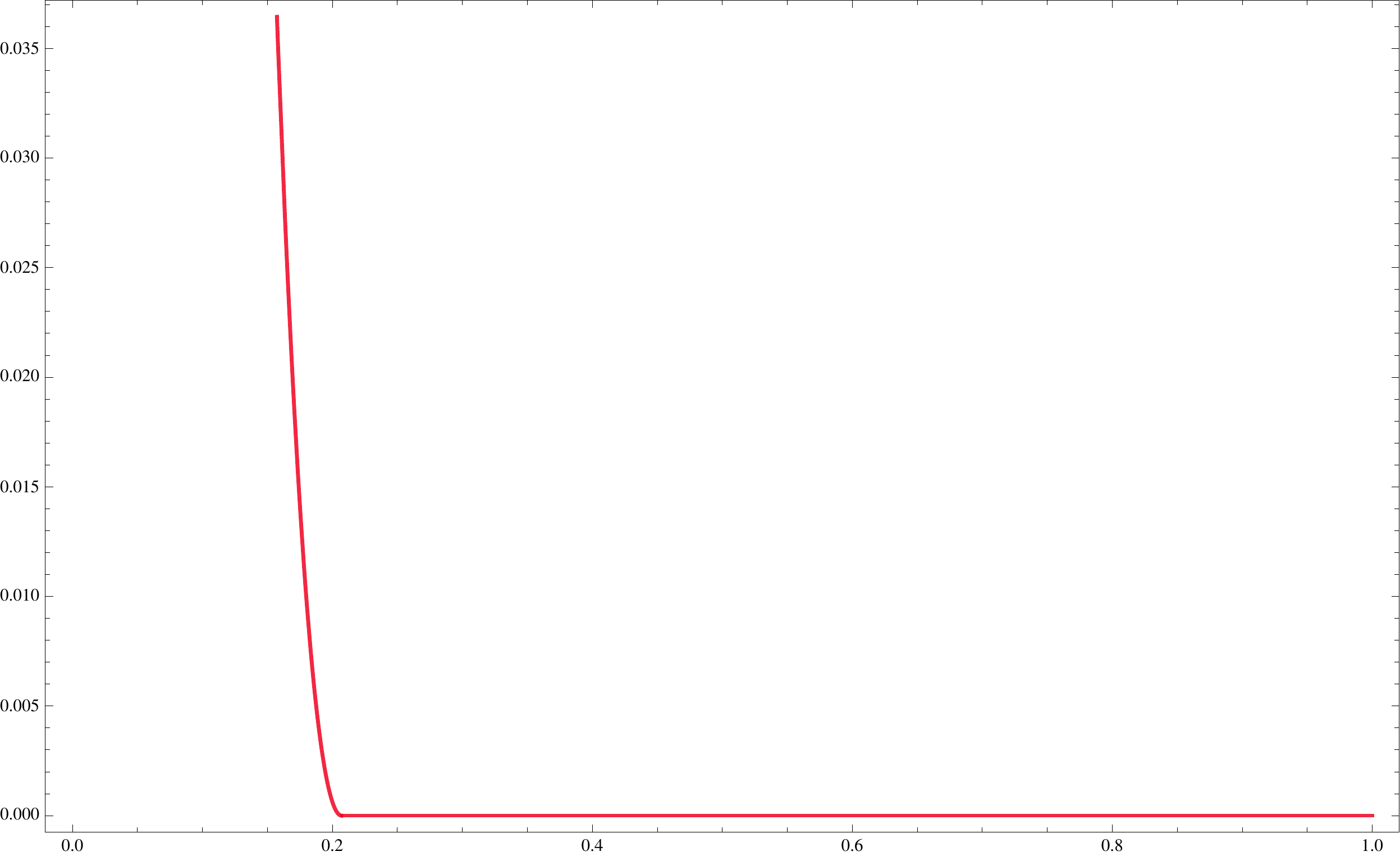}
\end{figure}

\begin{figure}[h!]
\caption{ $ h_{\ell_2} $ for $a=1 $ and $ \alpha=-12,\tau=288,\ell_2=1/2 $ }\label{fig3}
\includegraphics[scale=0.42]{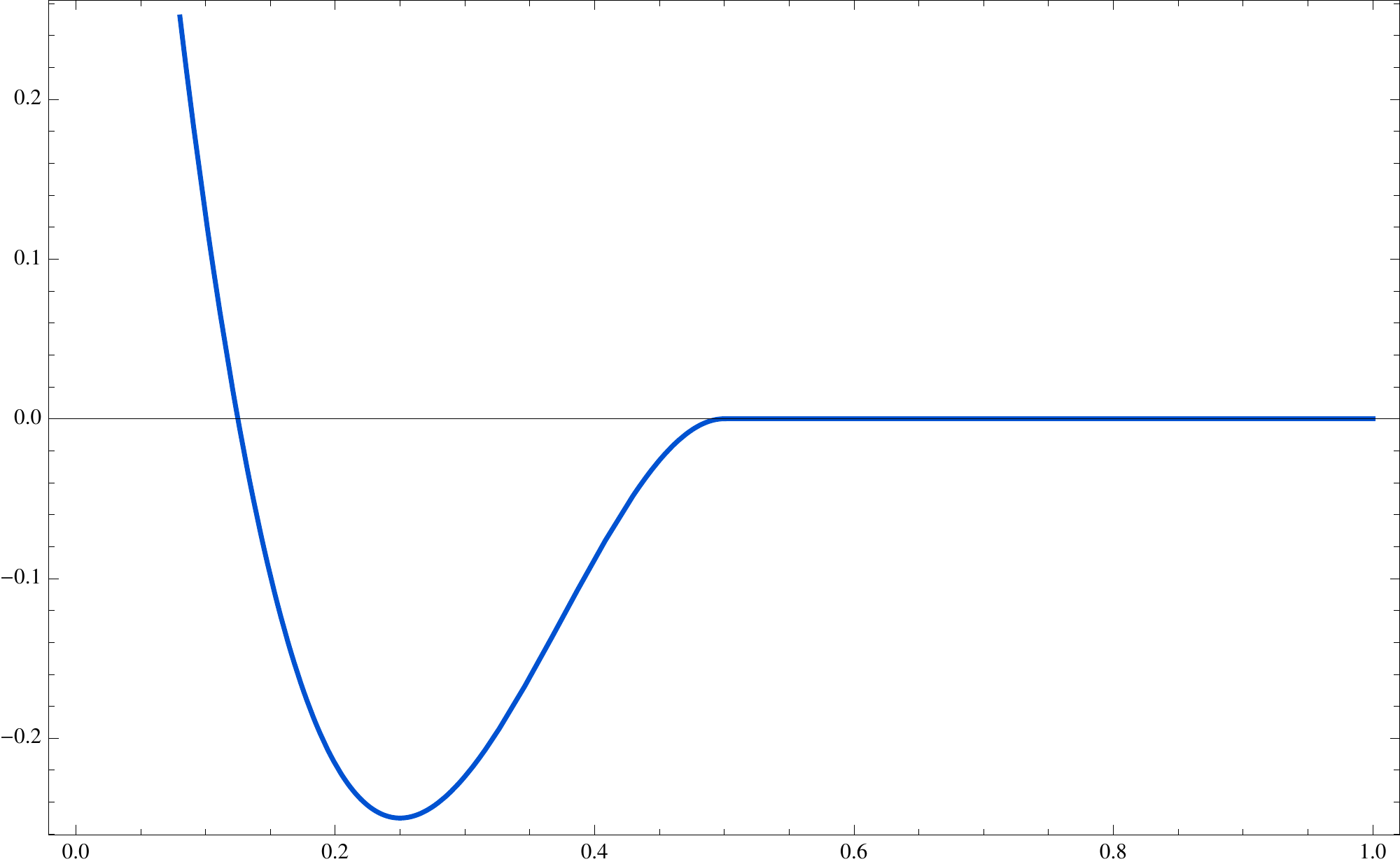}
\end{figure}

\begin{figure}[h!]
\caption{ $ h_{\ell_1} $ for $a=\alpha=1 $ and $\tau=288,\ell_1=6/100(1+\sqrt{101}) $ }\label{fig3a}
\includegraphics[scale=0.42]{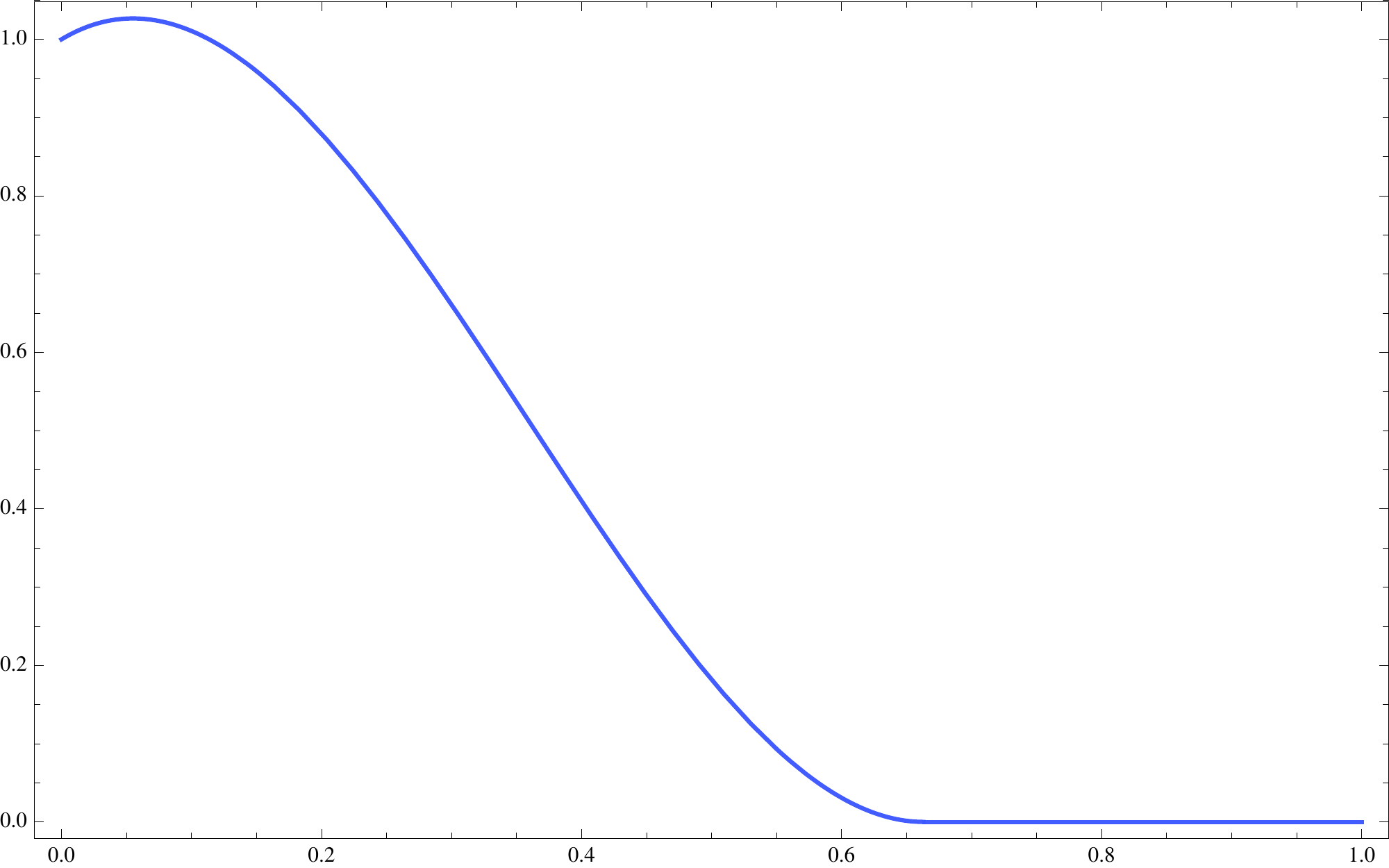}
\end{figure}

Note that $ \overline{h} $  does not pick up reward for any boundary condition $ \ba\not= \zero $ whereas  for $ \ba=\zero $ it takes the maximal reward. The function $ h_\ell $ picks up the reward in $ [\ell,1] $, see Figure~\ref{fig2},\ref{fig3},\ref{fig3a}. This motivates the following definitions. For any $ \tau\in\R $ and $ \ba\in\R^2 $ we let 
\begin{equation}\label{energyfct}
\Escr_{(\ba,\zero)}^\tau(\ell)=\Escr(h_{\ba,\zero}^{*,(0,\ell)})+\tau\ell, 
\end{equation}
and observe that for $ \tau=\tau(\eps) $
\begin{equation}
\Sigma^\eps_f(h_\ell)=\Escr^{\tau(\eps)}_{(\ba,\zero)}(\ell)-\tau(\eps).\end{equation}

Henceforth minimiser of $ \Sigma_f^\eps $ are given by functions of type  $ h_\ell $ only  if  $ \ell $ is a minimiser of the function $ \Escr_{(\ba,\zero)}^\tau $ in $ [0,1] $. We collect an analysis of the latter function in the next Proposition.

%

\begin{proposition}[\textbf{Minimiser for $ \boldsymbol{\Escr^\tau_{(\ba,\zero)}} $}]\label{minimafree}
\begin{enumerate}
\item[(a)] For $ \tau=0 $ the function  $ \Escr_{(\ba,\zero)}^0 $ is strictly decreasing with  $\lim_{\ell\to\infty}\Escr_{(\ba,\zero)}^0(\ell)=0$. 

\item[(b)]  For $ \tau>0 $ the function $ \Escr_{(a,0,\zero)}^\tau, a\not= 0 $, has one local minimum at $ \ell=\ell_1(\tau,a,0)=\sqrt{|a|}(18/\tau)^{1/4}$, and the function $ \Escr_{(0,\alpha,\zero)}^\tau ,\alpha \not=  0$, has one local minimum at $ \ell=\ell_1(\tau,0,\alpha)=\sqrt{\frac{2}{\tau}}|\alpha|$. In both cases there exist $ \tau_1(\ba) $ such that $ \ell_1(\tau,\ba) \le 1 $ for all $ \tau\ge\tau_1(\ba) $.


\item[(c)] For $ \tau>0 $ and $  \ba=(a,\alpha)\in\R^2 $ with $ w=|a|/|\alpha|\in (0,\infty) $ and $ s=\sign(a\alpha) $ the function $ \Escr^\tau_{(\ba,\zero)} $ has one local minimum at $ \ell=\ell_1(\tau,\ba)=\frac{1}{\sqrt{2\tau}}\big(|\alpha|+\sqrt{\alpha^2+6|a|\sqrt{2\tau}}\big)$ when $ s=1 $ , whereas for $ s=-1 $ the function $ \Escr^\tau_{(\ba,\zero)} $ has two local minima at $ \ell=\ell_1(\tau,\ba)= \frac{1}{\sqrt{2\tau}}\big(-|\alpha|+\sqrt{\alpha^2+6|a|\sqrt{2\tau}}\big)$ and $ \ell=\ell_2(\tau,\ba)=\frac{1}{\sqrt{2\tau}}\big(|\alpha|+\sqrt{\alpha^2-6|a|\sqrt{2\tau}}\big) $, where $ \ell_2 $ is a local minimum only if $ \tau\le\frac{\alpha^4}{72a^2} $. In all cases there are $ \tau_i(\ba) $ such that $ \ell_i(\tau,\ba)\le 1 $ for all $ \tau\ge \tau_i(\ba) $, $i=1,2$.
\end{enumerate}
\end{proposition}

From now we use the notation for $ \ell_1 $ and $ \ell_2 $ for the points where the functions $ h_{\ell_i},i=1,2$, pick up reward.
We shall study the zero sets of all minimiser, that is we need  to check if $ h_\ell $ has zeroes in $ [0,\ell) $ before picking up the reward in $ [\ell,1] $. 
\begin{lemma}\label{L:scaling}
Let $ a>0 $, then the functions $ h^{*,(0,\ell)}_{(a,\alpha,\zero)} $ with $ \alpha>0$, $h^{*,(0,\ell)}_{(0,\alpha,\zero)} $ with $\alpha\not=0 $, and $ h^{*,(0,\ell)}_{(a,-\alpha,\zero)} $ with $ \alpha\ell/a\in[0,3) $ have no zeroes in $ (0,\ell) $, whereas the functions $ h^{*,(0,\ell)}_{(a,-\alpha,\zero)} $ with $ \alpha\ell/a>3 $ have exactly one zero in $ (0,\ell) $. Analogous statements hold for $ a<0 $.
\end{lemma}
There is a qualitative difference between the minimiser $ h_{\ell_1} $ and $ h_{\ell_2} $ as the latter one has a zero before picking up the reward on $ [\ell_2,1] $, see Figure~\ref{fig3}.

%

In the following we write $ \eps_i(\ba) $ for the value of the reward with $ \tau(\eps_i(\ba))=\tau_i(\ba) $ such that $ \ell_i(\tau_i(\ba),\ba)\le 1, i=1,2 $.

\begin{theorem}[\textbf{Minimiser for $ \Sigma_f^\eps$}]\label{T:freemin}  
\begin{enumerate}
\item[(a)] If $ \ba=(a,0),a\not=0 $ or $ \ba=(0,\alpha), \alpha\not=0 $ or  $ w=|a|/|\alpha|\in (0,\infty) $ with $ s=\sign (a\alpha) =1 $, there exists $ \eps^*(\ba) >\eps_1(\ba) $ such that 
$$ 
\Mscr_f^\eps=\begin{cases} \{\overline{h}\} &, \mbox{ for } \eps<\eps^*(\ba),\\
\{\overline{h},h_{\ell_1}\}  \mbox{ with } \Sigma_f^{\eps^*}(\overline{h})= \Sigma_f^{\eps^*}(h_{\ell_1}) &, \mbox{ for }\eps= \eps^*(\ba),\\
\{h_{\ell_1} \}   &, \mbox{ for } \eps>\eps^*(\ba).\end{cases}
$$

%

\item[(b)] Assume $ w=|a|/|\alpha|\in (0,\infty) $ and $ s=\sign (a\alpha) = -1 $. There are $ \tau_0(\ba) > 0 $ and $ \tau_1^*(\ba) >0 $ and $ \tau^*_2(\ba)>0 $  such that the following statements hold.

 \begin{itemize}
\item[(i)]   Let $ \ba\in D_1:=\{\ba\in\R^2\colon w\in(0,\infty) \mbox{ and } \tau_0(\ba)>\tau^*_1(\ba)\} $.Then there exist $ \eps^*_{1,2}(\ba) >0 $ and $ \eps^*_2(\ba) >\eps_2(\ba)  $ with  $ \eps_2^*(\ba)<\eps^*_{1,2}(\ba) $ and $ \tau(\eps_{1,2}^*(\ba))=\tau_0(\ba) $ and $ \tau(\eps_2^*(\ba))=\tau_2^*(\ba) $ such that
$$
\Mscr_f^\eps=\begin{cases} \{\overline{h}\} &,\mbox{ for } \eps<\eps^*_2(\ba),\\
\{\overline{h},h_{\ell_2}\} \mbox{ with } \Sigma_f^{\eps^*_2(\ba)}(\overline{h})=\Sigma_f^{\eps^*_2(\ba)}(h_{\ell_2}) &, \mbox{ for }  \eps=\eps^*_2(\ba),\\
\{h_{\ell_2}\} &.\mbox{ for } \eps\in (\eps^*_2(\ba),\eps^*_{1,2}(\ba)),\\
\{h_{\ell_1},h_{\ell_2}\} \mbox{ with } \Sigma_f^{\eps^*_{1,2}(\ba)}(h_{\ell_1})=\Sigma_f^{\eps^{*}_{1,2}(\ba)}(h_{\ell_2}) &, \mbox{ for } \eps=\eps^*_{1,2}(\ba),\\
\{h_{\ell_1}\}   &, \mbox{ for } \eps>\eps^*_{1,2}(\ba).
\end{cases}
$$
\item[(ii)]  Let $ \ba\in D_2:=\{\ba\in\R^2\colon w\in (0,\infty) \mbox{ and } \tau_0(\ba)=\tau^*_1(\ba)\} $.  Then for $ \eps_{\rm c}^*(\ba) >\tau_2(\ba) $ with $ \tau(\eps^*_{\rm c}(\ba))=\tau_0(\ba) $,
$$
\Mscr_f^\eps=\begin{cases} 
\{\overline{h}\} &,\mbox{ for } \eps<\eps^*_{\rm c}(\ba),\\
\{\overline{h},h_{\ell_1},h_{\ell_2}\} \mbox{ with } \Sigma_f^{\eps^*_{\rm c}(\ba)}(\overline{h})=\Sigma_f^{\eps^*_{\rm c}(\ba)}(h_{\ell_1})=\Sigma_f^{\eps^*_{\rm c}(\ba)}(h_{\ell_2}) &, \mbox{ for } \eps=\eps^*_{\rm c}(\ba),\\
\{h_{\ell_1}\}   &, \mbox{ for } \eps>\eps^*_{\rm c}(\ba).
\end{cases}
$$

\item[(iii)]  Let $ \ba\in D_3:=\{\ba\in\R^2\colon w\in (0,\infty) \mbox{ and } \tau_0(\ba)<\tau^*_1(\ba)\} $. Then for $ \eps^*_1(\ba)>0 $ with $ \tau(\eps^*_1(\ba))=\tau^*_1(\ba) $,
$$
\Mscr_f^\eps=\begin{cases} \{\overline{h}\} &,\mbox{ for } \eps<\eps^*_1(\ba),\\
\{\overline{h},h_{\ell_1}\} \mbox{ with } \Sigma_f^{\eps^*_1(\ba)}(\overline{h})=\Sigma_f^{\eps^*_1(\ba)}(h_{\ell_1}) &, \mbox{ for } \eps=\eps^*_1(\ba),\\
\{h_{\ell_1}\}   &, \mbox{ for } \eps>\eps^*_1(\ba).\end{cases}
$$
\end{itemize}
\end{enumerate}
\end{theorem}

\begin{remark}
We have seen that the rate function $ \Sigma_f^\eps $ can have up to three distinct global minimisers.  See Figure~\ref{fig2}-\ref{fig3} for examples of these functions.  The minimiser in Figure~\ref{fig2} has no isolated zero before picking up the reward. Note that the existence of the minimiser (see Figure~\ref{fig3}) with a single zero before picking up a reward depends on the choice boundary conditions. This minimiser only exist if the gradient at $0$ has opposite sign of the value at zero. See Figure~\ref{fig3a} for an example when the gradient has the same sign as the value of the function at zero.  The minimiser $ h_{\ell_1} $ is the global minimiser if the reward is sufficiently large.
\end{remark}

\subsection{Dirichlet boundary}\label{Sec-MinD}
We consider Dirichlet boundary conditions on both sides given by the vector $ \br=(a,\alpha,b,\beta)=(\ba,\bb) $. In a similar way to Section~\ref{Sec-Minfree} for free boundary conditions on the right hand side  
we define functions $ h_{\ell,r}\in H_{\br}^2 $  for any $ \ell,r\ge 0 $ with $ \ell\le r $ and $ \ell+r\le 1 $ by
\begin{equation}
h_{\ell,r}(t)=\begin{cases} h_{\ba,\zero}^{*,(0,\ell)}(t) &, t\in [0,\ell),\\
0 &, t\in[\ell,1-r],\\
h_{\zero,\bb}^{*,(1-r,1)}(t) &, t\in(1-r,1].\end{cases}
\end{equation}
 Furthermore, we define the following energy function depending only on $ \ell $ and $ r $,
\begin{equation}
\begin{aligned}
E(\ell,r)&=\Escr(h_{\ba,\zero}^{*,(0,\ell)})+\Escr(h_{\zero,\bb}^{*,(1-r,1)})-\tau(\eps)(1-\ell-r),\\
\end{aligned}
\end{equation}
and using \eqref{energyfct} we get
\begin{equation}\label{energyfcts}
E(\ell,r)=\Escr_{(\ba,\zero)}^{\tau(\eps)}(\ell)+\Escr_{(\zero,b,\beta)}^{\tau(\eps)}(r)-\tau(\eps)= \Escr_{(\ba,\zero)}^{\tau(\eps)}(\ell)+\Escr_{(b,-\beta,\zero)}^{\tau(\eps)}(r)-\tau(\eps)  ,
\end{equation}
where $ \beta $ is replaced by $ -\beta $ due to symmetry, that is, using that $ h^{*,(1-r,1)}_{(\zero,\bb)}(t)=h^{*,(1-r,1)}_{(b,-\beta,\zero)}(2-r-t)=h^{*,(0,r)}_{(b,-\beta,\zero)}(1-t) $ for $ t\in[1-r,1]$. Hence
$$
\Sigma^\eps(h_{\ell,r})=E(\ell,r).  
$$
  

 \begin{figure}[h!]
\caption{$h_{\ell_1} $ for $a=b=1 $ and $ \alpha=-12,\beta=12,\tau=288,\ell_1=1/2(\sqrt{2}-1) $} \label{fig5}
\includegraphics[scale=0.42]{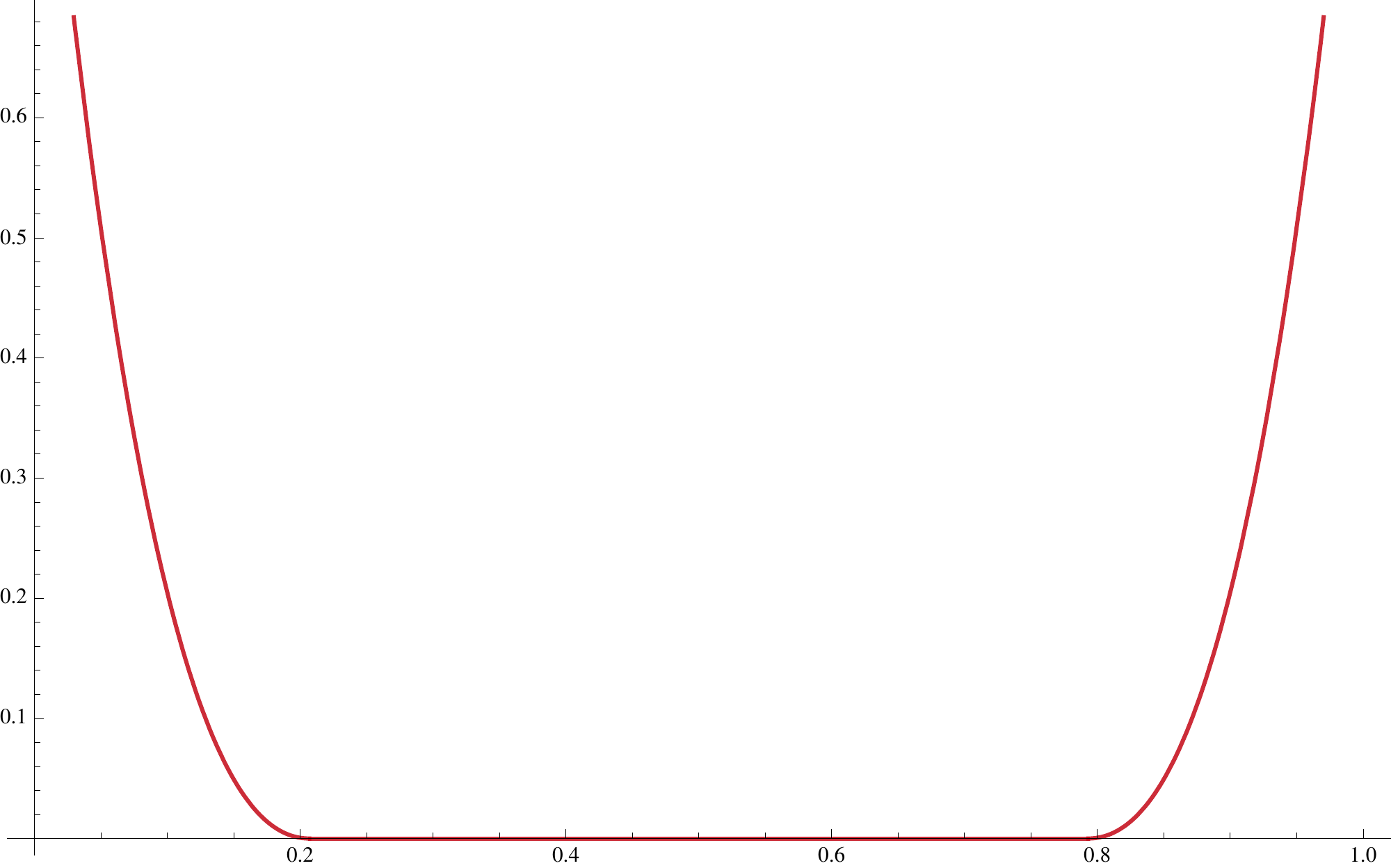}
\end{figure}

 \begin{figure}[h!]
\caption{ $ h_{\ell_2} $ for $a=b=1 $ and $ \alpha=-12,\beta=12,\tau=288,\ell_2=1/2 $ }\label{fig6}
\includegraphics[scale=0.42]{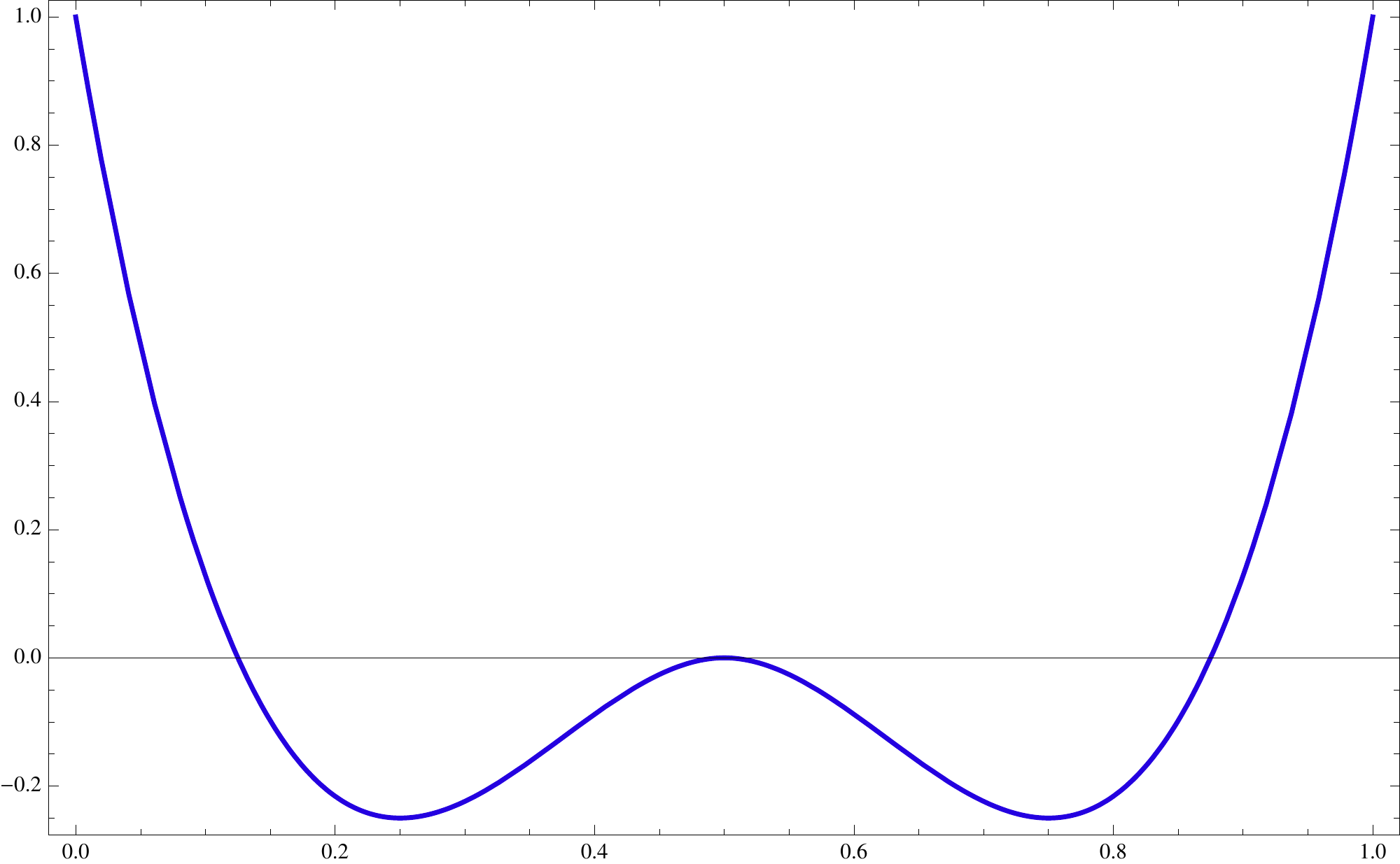}
\end{figure}

For given boundary $ \br=(a,\alpha,b,\beta) $ the function $ h^*_{\br} \in H_{\br}^2 $  given in Proposition~\ref{uniquemin} 
does not pick up any reward in $ [0,1] $.  

\begin{proposition}\label{minmainD}
For any Dirichlet boundary condition $ \br\in\R^4 $ the set $ \Mscr^\eps $ of minimiser of the rate function $ \Sigma^\eps $  in $ H_{\br}^2 $  is a subset of 
$$
\{h_{\ell,r}, h^*_{\br}\colon \ell+r\le 1\},
$$
where $ \ell $ and $ r $ are minimiser of $ \Escr^{\tau(\eps)} $ in Proposition~\ref{minimafree}.
\end{proposition}
Proposition~\ref{minmainD}  allows to reduce the optimisation of the rate function $ \Sigma^\eps $ to the minimisation of  the function $ E $ defined in \eqref{energyfcts} for $ 0\le \ell+r\le 1 $.  The general problem involves up to five parameters including  the boundary conditions $ \br\in\R^4 $ and the pinning free energy $ \tau(\eps) $ for the reward $ \eps$. It involves studying several different sub cases and in order to demonstrate  the key features of the whole  minimisation problem we study only  a special case in the following and  only outline how any general case can be approached. 

\noindent The \textit{symmetric case} $ \br=(a,\alpha,a,-\alpha) $:  It is straightforward to see that 
\begin{equation}\label{sym1}
\Sigma^\eps(h_{\ell_i,\ell_j})=E(\ell_i,\ell_j)=\Escr^{\tau(\eps)}_{(a,\alpha,\zero)}(\ell_i) + \Escr^{\tau(\eps)}_{(a,\alpha,\zero)}(\ell_j)-\tau(\eps), \quad i,j=1,2.
\end{equation}
Clearly the unique minimiser  $ h_{\br}^*(t)=a+\alpha t-\alpha t^2  $  of $ \Escr $  has the symmetry $ h^*_{\br}(1/2-t)=h^*_{\br}(1/2+t) $ for $ t\in[0,1/2] $. The function $ E  $ is not convex and thus we distinguish two different sets of parameter $ (a,\tau(\eps)) \in\R^3 $ according to whether (i) $ \ell_i(\tau(\eps),\ba)\le 1/2 $ for $ i=1,2 $;  or whether (ii) $ \ell_2(\tau(\eps),\ba)>1/2>\ell_1(\tau(\eps),\ba) $. There are no other cases for the parameter due to the condition $ \ell_1+\ell_2\le 1 $ and the fact that $ \ell_2(\tau(\eps),\ba)>\ell_1(\tau(\eps),\ba) $.

\noindent Parameter regime (i): 
$$
\Dscr_1:=\{(\ba,\tau)\in\R^3\colon \ell_1(\tau,\ba)\le 1/2 \,\wedge\,  \ell_2(\tau,\ba)\le 1/2 \mbox{ if }  \ell_2(\tau,\ba) \mbox{ is local minimum of } \Escr^\tau_{(\ba,\zero)}\}. 
$$ 
\noindent Parameter regime (ii): 
$$ 
\Dscr_2:=\{(\ba,\tau)\in\R^3\colon 1\ge \ell_2(\tau(\eps),\ba)>1/2>\ell_1(\tau(\eps),\ba)>0, \tau\le \frac{\alpha^4}{72a^2}\}. 
$$  

\noindent We shall define the following values before stating our  results.\\
There  are $ \eps_i(\ba) $ such that $ \ell_i(\tau(\eps),\ba)\le 1/2 $ for  all $ \eps\ge\eps_i(\ba) , i=1,2$.  We denote by $ \tau^*_1(\ba)=\tau(\eps^*_1(\ba)) $ the unique value of $ \tau $ such  that
\begin{equation}\label{sym2} 
\Escr^{\tau^*_1(\ba)}(\ell_1(\tau^*_1(\ba),\ba))-1/2\tau^*_1(\ba)=1/2 \Escr(h^*_{\br}).
\end{equation}
Likewise, we denote  $ \tau^*_2(\ba) $ the unique value of $ \tau $ such that  $  \Escr^{\tau^*_2(\ba)}(\ell_1(\tau^*_2(\ba),\ba))-1/2\tau^*_2(\ba)=1/2 \Escr(h^*_{\br}) $  when such a value exists in $ \R $ otherwise we put $ \tau_2^*(\ba)=\infty $. We denote  $ \tau_0(\ba) $ the unique zero  in Lemma~\ref{keylemma} (a) of the difference $ \Delta(\tau)=\Escr^{\tau}(\ell_1(\tau,\ba))- \Escr^{\tau}(\ell_2(\tau,\ba)) $. \\

%
%

\begin{theorem}[\textbf{Minimiser for $ \Sigma^\eps $, symmetric case}]\label{T:minD} Let $ \br=(a,\alpha,a,-\alpha) $. 
\begin{enumerate}
 \item[(a)] If $ \ba=(a,0),a\not=0, $ or $ \ba=(0,\alpha),\alpha\not=0,$ or $ w=|a|/|\alpha|\in(0,\infty) $ with $ \sign(a\alpha)=1 $ and $ \eps\ge \eps_1(\ba) $, then $ (\ba,\tau(\eps))\in\Dscr_1 $ and  there is $ \eps^*_1(\ba)>\eps_1(\ba) $ such that  
 $$
 \Mscr^\eps=\begin{cases} \{h^*_{\br}\} &, \mbox{ for } \eps <\eps^*_1(\ba),\\
 \{h^*_{\br}, h_{\ell_1,\ell_1}\} \mbox{ with } \Sigma^{\eps^*_1(\ba)}(h^*_{\br})=\Sigma^{\eps^*_1(\ba)}(h_{\ell_1,\ell_1}) &, \mbox{ for } \eps=\eps^*_1(\ba),\\
 \{h_{\ell_1,\ell_1}\} &, \mbox{ for } \eps>\eps^*_1(\ba).
 \end{cases}
 $$
%
%
%

 \item[(b)] Assume $ w=|a|/|\alpha|\in (0,\infty) $ and $ s=\sign (a\alpha) = -1 $. There are $ \tau_0(\ba) > 0 $ and $ \tau_1^*(\ba) >0 $ and $ \tau^*_2(\ba)>0 $ such that the following statements hold.

 \begin{itemize}
\item[(i)]   Let $ \ba\in D_1:=\{\ba\in\R^2\colon w\in(0,\infty) \mbox{ and } \tau_0(\ba)>\tau^*_1(\ba)\} $.Then there exists $ \widetilde{\eps}_{1,2}(\ba)>0 $ such that $ (\ba,\tau(\eps))\in \Dscr_2 $ for all $ \eps\in(\widetilde{\eps}_{1,2}(\ba),\eps_2(\ba)) $ and $ (a,\tau(\eps))\in\Dscr_1 $ for $ \eps\ge\eps_2(\ba) $. Then there exist $ \eps^*_{1,2}(\ba) >0 $ and $ \eps^*_2(\ba) >0 $ with  $ \eps_2^*(\ba)<\eps^*_{1,2}(\ba) $ and $ \tau(\eps_{1,2}^*(\ba))=\tau_0(\ba) $ and $ \tau(\eps_2^*(\ba))=\tau_2^*(\ba) $ such that
$$
\Mscr^\eps=\begin{cases} \{h^*_{\br}\} &,\mbox{ for } \eps<\widetilde{\eps}_{1,2}(\ba),\\
\{h^*_{\br},h_{\ell_2,\ell_1}\}\mbox{ with } \Sigma^{\eps}(h^*_{\br})\le \Sigma^{\eps}(h_{\ell_2,\ell_1}) \vee \Sigma^\eps(h^*_{\br})>\Sigma^\eps(h_{\ell_2,\ell_1}) &,\mbox{ for } \eps\in (\widetilde{\eps}_{1,2}(\ba),\eps_2(\ba)),\\ 
\{h^*_{\br},h_{\ell_2,\ell_2}\} \mbox{ with } \Sigma^{\eps^*_2(\ba)}(h^*_{\br})=\Sigma^{\eps^*_2(\ba)}(h_{\ell_2,\ell_2}) &, \mbox{ for }  \eps=\eps^*_2(\ba),\\
\{h_{\ell_2,\ell_2}\} &.\mbox{ for } \eps\in (\eps^*_2(\ba),\eps^*_{1,2}(\ba)),\\
\{h_{\ell_1,\ell_1},h_{\ell_2,\ell_2}\} \mbox{ with }  \Sigma^{\eps^*_{1,2}(\ba)}(h_{\ell_1,\ell_1})=\Sigma^{\eps^{*}_{1,2}(\ba)}(h_{\ell_2,\ell_2}) &, \mbox{ for } \eps=\eps^*_{1,2}(\ba),\\
\{h_{\ell_1,\ell_1}\}   &, \mbox{ for } \eps>\eps^*_{1,2}(\ba).
\end{cases}
$$
\item[(ii)]  Let $ \ba\in D_2:=\{\ba\in\R^2\colon w\in (0,\infty) \mbox{ and } \tau_0(\ba)=\tau^*_1(\ba)\} $.  Then there exists $ \widetilde{\eps}_{1,2}(\ba)>0 $ such that $ (\ba,\tau(\eps))\in \Dscr_2 $ for all $ \eps\in(\widetilde{\eps}_{1,2}(\ba),\eps_2(\ba)) $ and $ (a,\tau(\eps))\in\Dscr_1 $ for $ \eps\ge\eps_2(\ba) $.
Then there exists  $ \eps_{\rm c}^*(\ba) >0 $ with $ \tau(\eps^*_{\rm c}(\ba))=\tau_0(\ba) $ and $ \eps^*_{\rm c}(\ba)\ge\eps_2(\ba) $ such that 
$$
\Mscr^\eps=\begin{cases} 
\{h^*_{\br}\} &,\mbox{ for } \eps<\widetilde{\eps}_{1,2}(\ba),\\
\{h^*_{\br},h_{\ell_2,\ell_1}\}\mbox{ with } \Sigma^{\eps}(h^*_{\br})\le \Sigma^{\eps}(h_{\ell_2,\ell_1}) \vee \Sigma^\eps(h^*_{\br})>\Sigma^\eps(h_{\ell_2,\ell_1}) &,\mbox{ for } \eps\in (\widetilde{\eps}_{1,2}(\ba),\eps_2(\ba)),\\
\{h^*_{\br}\} &, \mbox{ for } \eps\in (\eps_2(\ba),\eps^*_{\rm c}(\ba)),\\
\{h_{\br}^*,h_{\ell_1,\ell_1},h_{\ell_2,\ell_2},h_{\ell_1,\ell_2},h_{\ell_2,\ell_1}\} \mbox{ with } \Sigma^{\eps^*_{\rm c}(\ba)}(h^*_{\br})= \Sigma^{\eps^*_{\rm c}(\ba)}(h_{\ell_1,\ell_1})= &\\ \Sigma^{\eps^*_{\rm c}(\ba)}(h_{\ell_2,\ell_2})=\Sigma^{\eps^*_{\rm c}(\ba)}(h_{\ell_1,\ell_2})=\Sigma^{\eps^*_{\rm c}(\ba)}(h_{\ell_2,\ell_1})  &, \mbox{ for } \eps=\eps^*_{\rm c}(\ba),\\
\{h_{\ell_1,\ell_1}\}   &, \mbox{ for } \eps>\eps^*_{\rm c}(\ba).
\end{cases}
$$

\item[(iii)]  Let $ \ba\in D_3:=\{\ba\in\R^2\colon w\in (0,\infty) \mbox{ and } \tau_0(\ba)<\tau^*_1(\ba)\} $. 
Then there exists $ \widetilde{\eps}_{1,2}(\ba)>0 $ such that $ (\ba,\tau(\eps))\in \Dscr_2 $ for all $ \eps\in(\widetilde{\eps}_{1,2}(\ba),\eps_1(\ba)) $ and $ (a,\tau(\eps))\in\Dscr_1 $ for $ \eps\ge\eps_1(\ba) $. Then there exists $ \eps^*_1(\ba)>\eps_1(\ba) $ with $ \tau(\eps^*_1(\ba))=\tau^*_1(\ba) $ such that
$$
\Mscr^\eps=\begin{cases} \{h^*_{\br}\} &,\mbox{ for } \eps<\widetilde{\eps}_{1,2}(\ba),\\
\{h^*_{\br},h_{\ell_2,\ell_1}\}\mbox{ with } \Sigma^{\eps}(h^*_{\br})\le \Sigma^{\eps}(h_{\ell_2,\ell_1}) \vee \Sigma^\eps(h^*_{\br})>\Sigma^\eps(h_{\ell_2,\ell_1}) &,\mbox{ for } \eps\in (\widetilde{\eps}_{1,2}(\ba),\eps_1^*(\ba)),\\
\{h^*_{\br},h_{\ell_1,\ell_1}\} \mbox{ with } \Sigma^{\eps^*_1(\ba)}(h^*_{\br})=\Sigma^{\eps^*_1(\ba)}(h_{\ell_1}) &, \mbox{ for } \eps=\eps^*_1(\ba),\\
\{h_{\ell_1,\ell_1}\}   &, \mbox{ for } \eps>\eps^*_1(\ba).\end{cases}
$$
\end{itemize}
  
 \end{enumerate}
 
 \end{theorem}

 \begin{remark}[\textbf{General boundary conditions}]
 For general boundary conditions $ \br=(a,\alpha,b,\beta) $ one can apply the same techniques as for the symmetric case. Thus minimiser of $ \Sigma^\eps $ are elements of 
 
 $$
 \{h^*_{\br}, h_{\ell,\ell},h_{r,r},h_{\ell,r},h_{r,\ell}\colon \ell+r\le 1\}.
 $$
 \end{remark}
\begin{remark}[\textbf{Concentration of measures}]
The large deviation principle in Theorem~\ref{THM-main} immediately implies the concentration properties for $ \gamma_N =\gamma_{N,\eps}^{\br} $ and $ \gamma_N= \gamma_{N,\eps}^{\ba} $:
\begin{equation}
\lim_{N\to\infty}\gamma_N(\dist_\infty(h_N,\Mscr^\eps)\le \delta)=1,
\end{equation}
for every $ \delta>0 $, where $ \Mscr^\eps=\{h^*\colon h^* \mbox{ minimiser of } I\} $ with $ I=I^\eps $ and $ I=I^\eps_f$, respectively, and $ \dist_\infty $ denotes the distance under $ \norm{\cdot}_\infty $. More precisely, for any $ \delta>0 $ there exists $ c(\delta)>0 $ such that
$$
\gamma_N(\dist_\infty(h_N,\Mscr^\eps)>\delta)\le \ex^{-c(\delta)N}
$$ for large enough $N$. We say that two function $ h_1,h_2\in\Mscr^\eps $ coexist in the limit $ N\to\infty $ under $ \gamma_N $ with probabilities $ \lambda_1,\lambda_2>0, \lambda_1+\lambda_2=1 $ when
$$
\lim_{N\to\infty}\gamma_N(\norm{h_N-h_i}_\infty\le \delta)=\lambda_i,\quad i=1,2,
$$ hold for small enough $ \delta>0 $. The same applies to the free boundary case on the right hand side and its set of minimiser $ \Mscr^\eps_f $. For gradient models with quadratic interaction (Gaussian) the authors in \cite{BFO} have investigated this concentration of measure problem and obtained statements depending on the dimension $ m $ of the underlying random walk (i.e. $ (1+m)$-dimensional models). The authors are using finer estimates than one employs for the large deviation principle, in particular the make  use of  a renewal property of the partition functions. In our setting of Laplacian interaction the renewal structure of the partition functions is different and requires different type of estimates. In addition, the concentration of measure problem requires to study all cases of possible minimiser. This is studied in (\cite{A16a}).
\end{remark}
\subsection{Proofs: Variational analysis}\label{Sec-minproof}
\subsubsection{Free boundary condition}
 
 \begin{proofsect}{Proof of Proposition~\ref{P:minfree}}
Suppose that $ h\in H_{\ba}^2 $ is not element of the set \eqref{superset}. It is easy to see that there is at least one function $ h^* $ in the set \eqref{superset} with
\begin{equation}\label{inequality}
\Sigma_f^\eps(h^*)<\Sigma_f^\eps(h).
\end{equation}
For $ \Sigma_f^\eps(h)<\infty $, we distinguish two cases. If $ |\Nscr_h|=0 $, then $ |\Nscr_{\overline{h}}|=0 $ and we get
$$
\Sigma_f^\eps(h)=\Escr(h)>0=\Escr(\overline{h})=\Sigma_f^\eps(\overline{h})
$$ by noting  that $ \overline{h} $ is the unique function with $ \Escr(\overline{h})=0 $. If $ |\Nscr_h|>0 $ we argue as follows. Let $ \ell $ be the infimum and $ r $ be the supremum of the accumulation points of $ \Nscr_h $, and note that $ \ell,r\in\Nscr_h $. Since $ |\Nscr_h\cap [\ell,r]^{\rm c}|=0 $ we have
$$
\begin{aligned}
\Sigma_f^\eps(h)=\Escr^{[0,\ell]}(h)+\Escr^{(\ell,r)}(h)-\tau(\eps)|\{t\in(\ell,r)\colon h(t)=0\}|+\Escr^{[r,1]}(h).
\end{aligned}
$$
As $ \ell,r\in\Nscr_h $ we have that  $ \dot{h}(\ell)=\dot{h}(r)=0$ as the differential quotient vanishes due  to the fact that $ \ell $ and $ r$ are accumulations points of $ \Nscr_h$. Thus the restrictions of $ h $ and $ h^*=h_\ell $ to $ [0,\ell] $ are elements of $ H_{(\ba,\zero)}^2 $. By the optimality of $ h_{(\ba,\zero)}^{*,(0,\ell)} $ inequality \eqref{inequality} is satisfied for $ h^*=h_\ell $. 
\qed
\end{proofsect}

\begin{proofsect}{Proof of Proposition~\ref{minimafree}}
The following scaling relations hold for $ \ell>0 $ (in our cases $ \ell\in(0,1)$) and $ \ba=(a,\alpha)$,
\begin{equation}
h^{*,(0,\ell)}_{(\ba,\zero)}(t)=h^{*,(0,1)}_{(a,\ell\alpha,\zero)}(t/\ell)\qquad\mbox{ for } t\in[0,\ell].
\end{equation}
Using this and Proposition~\ref{uniquemin} with $ \br=(a,\ell\alpha,0,0) $ we obtain
$$
\begin{aligned}
\Escr^{(0,\ell)}(h^{*,(0,\ell)}_{(\ba,\zero)})=\frac{1}{2\ell^3}\int_0^1\,\big(\ddot{h^{*}_{\br}}(t)\big)^2\,\d t=\frac{1}{\ell^3}\big(6a^2+6a\alpha\ell+2\alpha^2\ell^2\big),
\end{aligned}
$$ 
and thus
\begin{equation}\label{energies}
\begin{aligned}
\Escr^{\tau}_{(\ba,\zero)}(\ell)&=\frac{1}{\ell^3}\big(6a^2+6a\alpha\ell+2\alpha^2\ell^2\big)+\tau\ell,\\
\frac{\d}{\d\ell}\Escr^{\tau}_{(\ba,\zero)}(\ell)&=-\frac{18a^2}{\ell^4}-\frac{12a\alpha}{\ell^3}-\frac{2\alpha^2}{\ell^2}+\tau=-\frac{2}{\ell^4}(3a+\alpha\ell-\sqrt{\tau/2}\ell^2)(3a+\alpha\ell+\sqrt{\tau/2}),\\
\frac{\d^2}{\d\ell^2}\Escr^{\tau}_{(\ba,\zero)}(\ell)&=\frac{4}{\ell^5}(6a+\alpha\ell)(3a+\alpha\ell).
\end{aligned}
\end{equation}
The derivative has the following zeroes
$$
\begin{aligned}
\ell_{1/2}&=\frac{\alpha\pm\sqrt{\alpha^2+6a\sqrt{2\tau}}}{\sqrt{2\tau}},\\
\ell_{3/4}&=\frac{-\alpha\pm\sqrt{\alpha^2-6a\sqrt{2\tau}}}{\sqrt{2\tau}}.
\end{aligned}
$$
\noindent (a) Our calculations \eqref{energies} imply $ \frac{\d}{\d\ell} \Escr_{(\ba,\zero)}(\ell)< 0 $ for $ \ba\in\R^2\setminus\{\zero\} $ and $ \lim_{\ell\to\infty}\Escr_{(\ba,\zero)}(\ell)=0$. If $ a=\alpha=0 $, then $ \Escr_{(0,0,\zero)}(\ell)=0 $ for all $ \ell $.

 \noindent (b) If $ a\not=0 $ and $ \alpha=0 $ our calculations \eqref{energies} imply that the function has local minimum at $ \ell=\ell_1(\tau,\ba)=\sqrt{|a|}\big(\frac{18}{\tau}\big)^{1/4} $, whereas for $ a=0 $ and $ \alpha\not=0 $ the function has local minimum at $ \ell=\ell_1(\tau,\ba)=\sqrt{\tau/2}|\alpha|$.

\noindent (c) Let $ w=|a|/|\alpha|\in(0,\infty) $ and $ s=1$. Then \eqref{energies} shows that the function has a local minimum at $ \ell=\ell_1(\tau,\ba)=\frac{1}{\sqrt{2\tau}}(|\alpha|+\sqrt{\alpha^2+6|a|\sqrt{2\tau}}) $. If $ s=-1$ we get a local minimum at $ \ell=\ell_1(\tau,\ba)=\frac{1}{\sqrt{2\tau}}(-|\alpha|+\sqrt{\alpha^2+6|a|\sqrt{2\tau}}) $ and in case $ \tau\le\frac{\alpha^4}{72a^2} $ a second local minimum at $ \ell=\ell_2(\tau,\ba)=\frac{1}{\sqrt{2\tau}}(|\alpha|+\sqrt{\alpha^2-6|a|\sqrt{2\tau}}) $. Note that $ \ell_1(\tau,\ba)<\ell_2(\tau,\ba) $ whenever $ \ell_2(\tau,\ba) $ is local minimum. This follows immediately from the second derivative which is positive whenever $ \ell_i(\tau,\ba) \le \frac{3a}{|\alpha|} $ or $ \ell_i(\tau,\ba)\ge \frac{6a}{|\alpha|} $ for $ a>0>\alpha $ and $ i=1,2$.

\qed
\end{proofsect}

\begin{proofsect}{Proof of Lemma~\ref{L:scaling}}
We are using the scaling property
\begin{equation}\label{scaling}
\begin{aligned}
h^{*,(0,\ell)}_{(\ba,\zero)}(t)&=ah^{*,(0,1)}_{(1,sw^{-1}\ell,\zero)}(t/\ell)\quad \mbox{ for } a\not=0 \mbox{ and } t\in[0,\ell],s=\sign(a\alpha),w=|a|/|\alpha|,\\
h^{*,(0,\ell)}_{(\ba,\zero)}(t)&=h^{*,(0,1)}_{(0,\alpha\ell,\zero)}(t/\ell) \quad \mbox{ for } a=0 \mbox{ and } t\in[0,\ell],
\end{aligned}
\end{equation}
and show the following equivalent statements,
the functions $ h^{*}_{(1,\ell,\zero)} $ with $ \ell> 0$, $h^*_{(1,-\ell,\zero)} $ with $\ell\in\R\setminus\{0\} $, and $ h^*_{(0,\ell,\zero)} $ with $ \ell\in[0,3) $ have no zeroes in $ (0,1) $, whereas the functions $ h^*_{(1,-\ell,\zero)} $ with $ \ell>3 $ have exactly one zero in $ (0,1) $.
Thus we study the unique minimiser of $ \Escr $ given in Proposition~\ref{uniquemin}, that is, we consider first the functions $ h^*_{(1,s\ell,\zero)} $ for $ s\in\{-1,1\} $ and $ \ell>0 $. The function $ h^*_{(1,s\ell,\zero)} $ has a zero in $ (0,1) $ if and only if it has a local minimum at which it assumes a negative value. Its derivative has at most one zero in $ (0,1) $ as by Proposition~\ref{uniquemin} the derivative 
$$
\dot{h}^*_{\br}(t)=\ell s+2(-3-2\ell s)t+3(2+\ell s)t^2
$$ is  zero at $ t=1 $ for the boundary condition given by  $ \br=(1,s\ell,0,0) $. Now for $ s=1 $ the local extrema is a maximum as the function value at $ t=0 $ is greater than its value at $ t=1 $ and thus the derivative changes sign from positive to negative. For $ s=-1 $ and $ \ell\le 3 $ there is no local extrema as the first derivative is zero only at $ t=1 $ and has no second zero in $ (0,1) $ and the second derivative $ \ddot{h}^*_{(\br,\zero)}(1)=6-2\ell $   at $ t=1 $ is strictly positive. Thus the derivative takes only negative values in $ [0,1) $ and is zero at $t=0 $. For $ s=-1 $ and $ \ell>3 $ there is a local minimum as the second derivative at $ t=1 $ is now strictly negative implying that the first derivative changes sign from negative to positive and thus has a zero at which the function value is negative.  The functions $ h^*_{(0,\ell,\zero)} $ have no zero in $ (0,\ell) $ for $ \ell\not=0 $ by definition as the only zeroes are $ t=0 $ and $ t=1 $.
\qed
\end{proofsect}

\begin{proofsect}{Proof of Theorem~\ref{T:freemin}}
(a) (i) Let $ \alpha=0 $ and $ a\not=0 $. Note that $ \ell_1(\tau,\ba)\le 1 $ if and only if $ \tau\ge 18|a|^2 $. Let $ \eps_1(\ba) $ be the maximum of $ \eps_{\rm c} $ and this lower bound. We write $ \tau=\tau(\eps) $. Now $ \Sigma_f^\eps(h_{\ell_1})=0 $ if and only if $\tau=\tau^* $ with
$$
\frac{6\sqrt{|a|}}{18^{3/4}}+18^{1/4}\sqrt{|a|}=(\tau^*)^{1/4},
$$ and we easily see that $ \Sigma_f^\eps(h_{\ell_1(\tau,\ba)})<0 $ for all $ \tau>\tau^* $.

\noindent (ii) Now let $ a=0 $ and $ \alpha\not=0 $. Note that $ \ell_1(\tau,\ba)\le 1 $ if and only if $ \sqrt{\tau}\ge\sqrt{2}|\alpha| $ and thus let $ \eps_1(\ba) $ be the maximum of $ \eps_{\rm c} $ and $ 2|\alpha|^2 $. Now $ \Sigma_f^\eps(h_{\ell_1})=0 $ if and only if $\tau=\tau^*(\ba):=8|\alpha|^2 $, and
$ \Sigma_f^\eps(h_{\ell_1(\tau,\ba)})<0 $ as $ \d/\d\tau(\Sigma_f^\tau(h_{\ell_1(\tau,\ba)}) <0$ for all $ \tau>\tau^* $.

\noindent (iii) Now let $ s=\sign(a\alpha) =1$ and assume that  $ a,\alpha>0 $ (the case $ a,\alpha<0 $ follows analogously). As $ \ell_1(\tau,\ba) $ is decreasing in $ \tau>0 $ there is $ \eps_1(\ba)\ge \eps_{\rm c} $ such that $ \ell_1(\tau,\ba) \le 1 $ for all $ \tau\ge \tau_1(\ba) $. Lemma~\ref{keylemma} (b) shows that there exists $ \eps^*_1(\ba) $ such that $ \Sigma_f^{\eps_1^*(\ba)}(h_{\ell_1(\tau_1^*,\ba)})=0 $ and the uniqueness of that zero gives $ \Sigma_f^\eps(h_{\ell_1(\tau(\eps),\ba)})<0 $ for all $ \eps>\eps_1^*(\ba) $.

\noindent (b) Let $ s=\sign(a\alpha)=-1 $ and assume  $ a>0>\alpha $ (the other case follows analogously). Clearly  we have $ \Sigma_f^{\eps_i(\ba)}(h_{\ell_i(\tau(\eps_i(\ba)),\ba)})>0 $ as $ \ell_i(\tau(\eps_i(\ba)),\ba)=1 $ for $ i=1,2$, and for any $ \eps>\eps_i(\ba) $ we have $ \ell_i(\tau(\eps),\ba)<1 $ and thus
$$
\frac{\d}{\d\tau} f_i(\tau)=\ell_i(\tau,\ba)-1<0  \mbox{ where we write  } f_i(\tau)=\Sigma_f^\tau(h_{\ell_i(\tau,\ba)}), \quad i=1,2.
$$
Furthermore, due to Lemma~\ref{keylemma} there is a unique $ \tau_0=\tau_0(\ba) $ such that
$$
f_1(\tau)\ge f_2(\tau) \quad \mbox {for } \tau\le \tau_0 \mbox{ and } f_1(\tau) \le f_2(\tau) \mbox{ for } \tau\ge \tau_0 \mbox{ and } f_1(\tau_0)=f_2(\tau_0). 
$$
We thus know that $ f_1 $ is decreasing and that $ f_1(\tau)\to -\infty $ for $ \tau\to\infty $. As $ f_1(\tau_1(\ba))>0 $ there must be at least one zero which we denote $ \tau^*_1(\ba) $ which we write as $\tau(\eps_1^*(\ba)) $. The uniqueness of $ \tau_1^*(\ba) $ is shown in Lemma~\ref{keylemma} (b). Similarly, we denote by $ \tau_2^*(\ba) $ the zero of $ f_2 $ when this zero exists (otherwise we set it equal to infinity), and one can show uniqueness of this zero in the same way as done for $ \tau_1^*(\ba) $ in Lemma~\ref{keylemma} (b).  We can now distinguish three cases according to the sign of the functions $f_1 $ and $ f_2 $ at the unique zero $ \tau_0 $ of the difference $ \Delta=f_1-f_2 $. That is, we distinguish  whether $ \tau_0(\ba) $ is greater, equal or less the unique zero $ \tau^*_1(\ba) $ of $ f_1 $. 

\noindent (i) Let $ \ba\in D_1:=\{\ba\in\R^2\colon w\in(0,\infty) \mbox{ and } \tau_0(\ba)>\tau^*_1(\ba)\} $. Then  $ f_1(\tau_0(\ba))=f_2(\tau_0(\ba))<0 $ and thus  $ \tau^*_2(\ba) $ exists and satisfies $ \tau_2^*(\ba)<\tau^*_1(\ba) $. This implies immediately the statement by choosing $ \eps_{1,2}^*(\ba) $ and $ \eps_2^*(\ba) $ such that $ \tau(\eps^*_{1,2}(\ba))=\tau_0(\ba) $ and $ \tau(\eps^*_2(\ba))=\t_2^*(\ba) $.

\smallskip
\noindent (ii) Let $ \ba\in D_2:=\{\ba\in\R^2\colon w\in (0,\infty) \mbox{ and } \tau_0(\ba)=\tau^*_1(\ba)=\tau^*_2(\ba)\} $.  Then $ f_1(\tau_0(\ba))=f_2(\tau_0(\ba))=0 $ and thus for $ \eps^*_{\rm c}(\ba) $ with $ \tau(\eps^*_{\rm c}(\ba))=\tau_0(\ba) $ we get 
$ \Sigma_f^{\eps^*_{\rm c}(\ba)}(h_{\ell_1})=\Sigma_f^{\eps^*_{\rm c}(\ba)}(h_{\ell_2})=\Sigma^{\eps^*_{\rm c}(\ba)}_f(\overline{h})=0 $. Then Lemma~\ref{keylemma} (a) gives $ \Sigma_f^\eps(h_{\ell_1})<\Sigma_f^\eps(h_{\ell_2})<0 $ for all $ \eps>\eps^*_{\rm c}(\ba) $.

\smallskip
\noindent (iii) Let $ \ba\in D_3:=\{\ba\in\R^2\colon w\in (0,\infty) \mbox{ and } \tau_0(\ba)<\tau^*_1(\ba)\} $. Then $ f_1(\tau_0(\ba))=f_2(\tau_0(\ba))>0 $ and for $ \eps^*_1(\ba) $  with $ \tau(\eps^*_1(\ba))=\tau^*_1(\ba) $ we get $ \Sigma_f^{\eps^*_1(\ba)}(h_{\ell_1})=\Sigma_f^{\eps^*_1(\ba)}(\overline{h})=0 $ and $ \Sigma_f^\eps(h_{\ell_1})<0 $ and $ \Sigma^\eps_f(h_{\ell_1})<\Sigma_f^\eps(h_{\ell_2}) $ for $ \eps> \eps^*_1(a) $.

\qed
\end{proofsect}

\begin{lemma}\label{keylemma}
\begin{enumerate}
\item[(a)] 
 For any $ \ba\in\R^2 $ with $ w=|a|/|\alpha| \in(0,\infty)$  and $ \tau \in(0,\frac{\alpha^4}{72a^2}] $ the function
$$
\Delta(\tau):=\Escr^{\tau}(\ell_1(\tau,\ba))-\Escr^{\tau}(\ell_2(\tau,\ba))
$$
has a unique zero called $ \tau_0 $, is strictly decreasing and strictly positive for $ \tau<\tau_0 $.
\item[(b)] For any $ \ba\in\R^2$ with $ w\in (0,
\infty) $ there is a unique solution of 
\begin{equation}\label{zeroeq}
\Sigma_{f}^{\eps}(h_{\ell_1(\tau(\eps),\ba)})=0,\qquad  \tau=\tau(\eps)\ge\tau_1(\ba),
\end{equation} which we denote by $ \tau^*_1=\tau(\eps^*_1(\ba)) $.
\end{enumerate}

\end{lemma}
\begin{proofsect}{Proof of Lemma~\ref{keylemma}}
(a) The sign of the function $ \Delta $ is positive  for  $ \tau\to 0 $  whereas the sign is negative if $ \tau=\frac{\alpha^4}{72a^2} $. Hence, the continuous function $\Delta $ changes its sign and must have a zero. We obtain the uniqueness of this zero by showing that the function $ \Delta $  is strictly decreasing. For fixed $  \tau $ we have (Proposition~\ref{minimafree})
$$
\frac{\d}{\d\ell}\Escr^{\tau}(\ell)=0,\qquad \mbox{ for } \ell=\ell_2(\tau,\ba) \mbox{ or } \ell=\ell_2(\tau,\ba).
$$
The functions $ \Escr^{\tau}(\ell_i(\tau,\ba)) $ are rational functions of $ \ell_i(\tau,\ba) $ and depend explicitly on $ \tau $ as well. Thus the chain rule gives
$$
\frac{\d}{\d \tau} \Escr^{\tau}(\ell_i(\tau))=\ell_i(\tau).\qquad i=1,2 \mbox{ and } \tau\in(0,\frac{\alpha^4}{72a^2}].
$$ As $ \ell_1(\tau,\ba)<\ell_2(\tau,\ba) $ the first derivative of $ \Delta $ is negative on $ (0,\frac{\alpha^4}{72a^2}] $.

(b) We let $ \tau^*_1=\tau(\eps^*_1(\ba)) $ denote the solution of \eqref{zeroeq}. As the rate function is strictly positive for vanishing $ \tau $ and $ \lim_{\tau\to\infty}\Sigma_f^\tau(\ell_1(\tau,\ba))=-\infty $ we shall check whether there is a second solution to \eqref{zeroeq}. Suppose there are $ \tau(\eps)>\tau(\eps^\prime) $ solving  \eqref{zeroeq}  with
\begin{equation}\label{zeq1}
\Sigma_f^\eps(h_{\ell_1(\tau(\eps),\ba)})=\Sigma_f^{\eps^\prime}(h_{\ell_1(\tau(\eps^\prime),\ba)}).
\end{equation}
For fixed $ \ell $ the function $ \tau\mapsto\Sigma_f^\eps(h_\ell)=\Escr^\tau_{(\ba,\zero)}(\ell)-\tau  $ is strictly decreasing and thus 
\begin{equation}\label{zeq2}
\Sigma_f^{\eps^\prime}(h_\ell)>\Sigma_f^\eps(h_\ell)\qquad \mbox{ for } \ell=\ell_1(\tau(\eps^\prime)).
\end{equation}
Now Proposition~\ref{minimafree} gives
$$
\begin{aligned}
\Sigma_f^{\eps}(h_{\ell_1(\tau(\eps^\prime),\ba)})\ge \min_{\ell\in(0,1)}\Sigma_f^{\eps}(h_\ell)=\Sigma_f^\eps(h_{\ell_1(\tau(\eps),\ba)}).
\end{aligned}
$$ Combining \eqref{zeq1} and \eqref{zeq2} we arrive at a contradiction and thus the solution of \eqref{zeroeq} is unique. Hence
$ \Sigma_f^\eps(h_{\ell_1(\tau(\eps),\ba)}) <0 $ for all $ \tau(\eps)>\tau^*_1=\tau(\eps^*_1(\ba))$.

\qed
\end{proofsect}

\subsubsection{Dirichlet boundary conditions}

\begin{proofsect}{Proof of Proposition~\ref{minmainD}}
We argue as in our proof of Proposition~\ref{P:minfree} using \eqref{energyfcts} observing that for any $ h\in H_{\br}^2 $ with $ \ell $ being the infimum of accumulation points of $ \Nscr_h $ and $ 1-r $ being the corresponding supremum,
$$
\begin{aligned}
\Sigma^\eps(h)&=\Escr^{(0,\ell)}(h)-\tau(\eps)(1-\ell-r)+\Escr^{(1-r,1)}(h)=\Escr^{\tau(\eps)}_{(\ba,\zero)}(\ell) +\Escr^{\tau(\eps)}_{b,-\beta,\zero)}(r)-\tau(\eps),\\
&=E(\ell,r).
\end{aligned}
$$ 
The second statement follows from the Hessian of $ E $ being the product 
$$
\frac{\partial^2}{\partial \ell^2} \Escr^{\tau(\eps)}_{(\ba,\zero)}(\ell)\frac{\partial^2}{\partial r^2} \Escr^{\tau(\eps)}_{(b,-\beta,\zero)}(r)
$$
of the second derivatives of the functions $ \Escr^{\tau(\eps)} $ (see Proposition~\ref{minimafree}).
\qed
\end{proofsect}

\begin{proofsect}{Proof of Theorem~\ref{T:minD}}
(a): We first note that due to convexity of $ \Escr $ the solutions $ h^*_{\br} $ for boundary conditions $ \br=(a,\alpha,a,-\alpha) $ are symmetric with respect to the $1/2- $ vertical line.
Furthermore, in all three cases of (a) only $ \ell_1(\tau,\ba) $ is a minimiser of $\Escr^\tau $ and thus of $E$ due to symmetric boundary conditions and thus the Hessian (see above) of the energy function $ E $ \eqref{energyfcts} is positive implying convexity.  Henceforth, when $ \ell_1(\tau,\ba) $ is a minimiser of $ E$ the corresponding minimiser function (see Proposition~\ref{minmainD}) of the rate function has to be symmetric with respect to the $1/2-$ vertical line. These observations immediately give the proofs for all three cases in (a)  of Theorem~\ref{T:minD} because symmetric minimiser exist only if $ \ell_1(\tau,\ba)\le 1/2 $.  Hence we conclude with Theorem~\ref{T:freemin} and $ \ell_1(\tau,\ba)\le 1/2 $ for $ \eps\ge \eps_1(\ba) $ using the existence  of $\tau^*_1(\ba) $ solving \eqref{sym2}. The existence and uniqueness of $ \tau_1^*(\ba) $ can be shown using an adaptation of  Lemma~\ref{keylemma} (b). 

\noindent We are left to show all three sub cases (i)-(iii) of  (b) in Theorem~\ref{T:minD}.  In all these cases we argue differently depending on the parameter regime.  If  $ (\ba,\tau)\in\Dscr_1 $ we can argue as follows. If $ \ell_1(\tau,\ba) $ and $ \ell_2(\tau,\ba) $ both exist and  are minimiser of the energy function $ E$ we obtain convexity as above (the mixed derivatives vanish due to the fact that $ E$ is a sum of functions of the single variables). Then we can argue as above and conclude with our statements for all there sub cases for parameter regime $ \Dscr_1 $ with $ \ell_i(\tau,\ba)\le 1/2 , i=1,2 $.

The only other case for the minimiser $ \ell_1(\ba,\tau(\eps)) $ of $ \Escr^\tau $ is $ 1\ge \ell_2(\tau(\eps),\ba)>1/2>\ell_1(\tau(\eps),\ba)>0 $ which gives a candidate for minimiser of $ \Sigma^\eps $ which is not symmetric with respect to the $1/2-$ vertical line. It is clear that at the boundary of $ \Dscr_2 $, namely $ \ell_1+\ell_2=1 $, we get $ \Escr(h^*_{\br})<E(\ell_2(\tau(\eps),\ba),\ell_1(\tau(\eps),\ba))$. Depending on the values of the boundary conditions and the value of $ \tau(\eps) $ the minimiser can be either $ h^*_{\br} $ or the non-symmetric function $ h_{\ell_2,\ell_1} $, or both. As outlined in \cite{Antman} for elastic rods which pose similar variational problems there are no general statements about the minimiser in this regime, for any given values of the parameter one can check by computation which function has a lower numerical value.
\qed
\end{proofsect}

\section{Proofs: Large deviation principles}\label{Sec-PLDP}
In this chapter the  proofs for the large deviation theorems are presented. In Section~\ref{Sec-LDPnopinning} we prove the extension of Mogulskii's theorem to integrated random walks and integrated random walk bridges. 
%
In Section~\ref{Sec-LDPpinning} we prove the main large deviation result,  Theorem~\ref{THM-main}, for  models with  pinning. 
The proof of the lower LDP bound in Section~\ref{Sec-LDPpinning} relies on the Gaussian LDP via Lemma~\ref{le:GaussLDP}.  The proof of the upper LDP bound  
relies on a stronger Gaussian large deviation bound in the form of the Gaussian isoperimetric inequality presented in Lemma~\ref{le:Isoperi}.

\subsection{Sample path large deviation for integrated random walks and integrated random walk bridges}\label{Sec-LDPnopinning}
We  show Theorem~\ref{LDP-nopinning} by using the contraction principle and an adaptation of Mogulskii's theorem (\cite[Chapter 5.1]{DZ98}). 

\noindent (a) Recall the integrated random walk representation in Section~\ref{Sec-LDPnp} and define a family of random variables indexed by $t$ as
$$
\widetilde{Y}_N(t)=\frac{1}{N}Y_{\floor{Nt}+1},\quad 0\le t\le 1,
$$
and let $ \mu_N $ be the law of $ \widetilde{Y}_N $ in $ L_\infty([0,1]) $. From Mogulskii's theorem \cite[Theorem~5.1.2]{DZ98} we obtain that $ \mu_N $ satisfy in $ L_\infty([0,1]) $ the LDP with the good rate function
$$
I^M(h)=\begin{cases} \int_0^1\L^*(\dot{h}(t))\,\d t &, \mbox{ if } h\in\Ascr\Cscr, h(0)= \alpha,\\
\infty & \mbox{ otherwise},\end{cases}
$$
where $ \Ascr\Cscr $ denotes the space of absolutely continuous functions. 
The empirical profiles $ h_N $ are functions of the integrated random walk $ (Z_n)_{n\in\N_0} $ (see Proposition~\ref{P:CD}), and 
$$
\begin{aligned}
h_N(t)&=\frac{1}{N^2}Z_{\floor{Nt}}+\frac{1}{N^2}\int_{\frac{\floor{Nt}}{N}}^t \,\big(Z_{\floor{N
s}+1}-Z_{\floor{Ns}}\big)\,\d s=\frac{1}{N}\int_0^t\,\widetilde{Y}_N(s)\,\d s.
\end{aligned}
$$
The contraction principle applied to the integral mapping immediately immediately gives the LDP for the empirical profiles $h_N$.
The rate function for this LDP  is given as the following infimum
$$
J(h)=\inf_{g\in\Sscr_h} I^M(g), \quad \mbox{ with } \Sscr_h=\{g\in L_\infty([0,1])\colon \int_0^t\,g(s)\,\d s=h(t), t\in [0,1]\}.
$$
If either $ \dot{h}(0)\not= \alpha $ or $ h$ is not differentiable, then $ \Sscr_h=\emptyset $. In the other cases one obtains $ \Sscr_h=\{\dot{h}\} $, and therefore $ J\equiv \Escr_f $. This proves part  (a) of Theorem~\ref{LDP-nopinning}.

\noindent (b) In the Gaussian case the LDP  can be shown by Gaussian calculus (e.g., \cite{DS89}), or by employing the  contraction principle for the Gaussian integrated random walk bridge. The explicit distribution of the Gaussian bridge leads to the follows mapping. We only sketch this approach for illustrations. 
For simplicity choose the boundary condition $ \br=\zero$ and $ \ba=(0,0)$. The cases for non-vanishing boundary conditions follow analogously. Then
\begin{equation}\label{bridgemeasures}
\P^{\zero}=\P^{(0,0)}\circ B_N^{-1},
\end{equation}
where for $ Z=(Z_1,\ldots,Z_{N+1}) $,
$$
B_N(Z)(x)=Z_x-A_N(x,Z_N,Z_{N+1}-Z_N),\quad x\in\{1,2,\ldots,N+1\},
$$
and 
$$
A_N(x,u,v)=\frac{1}{N(N+1)(N+2)}\Big(x^3(-2u+vN)+x^2(3uN+vN-vN^2)+x((2+3N)u-N^2v\Big).
$$
Clearly, $ B_N(Z)(N)=B_N(Z)(N+1)=0 $. Now we see that the integrated random walk bridge distribution on the left hand side of \eqref{bridgemeasures}  is given by the integrated random distribution via  the continuous mapping $B_N$.  Therefore we can apply our reasoning in part (a) and another application of the contraction principle leads to the statement. Note that the explicit map $ B_N $ is only given for quadratic potentials, for more general potentials a different techniques will be required.

\subsection{Sample path large deviation for pinning models}\label{Sec-PLDPpinning}

In the following we will prove Theorem~\ref{THM-main} for the case of Dirichlet boundary conditions. We concentrate on the Dirichlet boundary case and only briefly comment on the (minor) difference in the case of free boundary conditions on the right towards the end of Sections~\ref{mainlower} and \ref{mainupper}. 
In Section~\ref{mainlower} we show the large deviation lower bound and in Section~\ref{mainupper} the corresponding upper bound.
It will be convenient to work in a slightly different normalisation. Instead of \eqref{LDPBounds} we will show that
 \begin{align}
 \liminf_{N\to\infty}\frac{1}{N}\log \frac{\ZNe(\br)}{Z_{N}(\zero)}\gNr(h_N\in \Oscr)&\ge -\inf_{h\in\Oscr} \Sigma(h),\label{e:LDPlowbou1} \\
 \limsup_{N\to\infty}\frac{1}{N}\log \frac{\ZNe(\br)}{Z_{N}(\zero)}\gNr(h_N\in\Kscr)&\le -\inf_{h\in\Kscr} \Sigma(h),\label{e:LDPuppbou2} 
 \end{align}
  where $\ZNe(\br)$ is the partition function introduced in \eqref{gibbsdist} with Dirichlet boundary condition given in \eqref{bc}  and $Z_{N}(\zero)$ is the partition function of the same model with pinning strength $\eps =0$ and Dirichlet boundary condition zero. Note for later use that exact formulae for the Gaussian partition function  $Z_{N}(\zero)$ are presented in Appendix~\ref{AppB}.   
  Once the bounds \eqref{e:LDPlowbou1} and \eqref{e:LDPuppbou2} are established, they can be applied to the full space $\Oscr =\Kscr= \Cscr([0,1];\R) $  implying 
 \begin{align*}
 \lim_{N\to\infty}\frac{1}{N}\log \frac{\ZNe(\br)}{Z_{N}(\zero)}&= -\inf_{h\in H} \Sigma(h),
 \end{align*}
so that \eqref{LDPBounds} follows.

\subsubsection{Proof of the lower bound in  Theorem~\ref{THM-main}}\label{mainlower} Fix  $g\in H_{\br}^2 $ and $\delta>0$. 
 We establish the lower bound \eqref{e:LDPlowbou1} in the form
\begin{align}\label{e:LDPlowbou} 
  \liminf_{N\to\infty}\frac{1}{N}\log \frac{\ZNe(\br)}{Z_{N}(\zero)}\gNr( \|h_N -g\|_\infty < \delta )&\ge - \Sigma(g).
 \end{align}
 
 \begin{proofsect}{Reduction to ``well behaved'' $g$} Recall that by Sobolev embedding any $g \in H_{\br}^2$ is automatically $\Cscr^1([0,1])$ with $\frac12$-H\"older continuous first derivative. We can write 
 \begin{align*}
 \{ t \in [0,1] \colon g(t) = 0\} = \bar{\Nscr} \cup \Nscr,
\end{align*} 
  where $\bar{\Nscr}$ is the set of \emph{isolated} zeros 
   \begin{align*}
 \bar{\Nscr} = \{ t \in [0,1] \colon g(t ) = 0 \, \text{ and } g \text{ has no further zeros in an open interval around }t \},
 \end{align*}
 and where $ \Nscr $ is the set of all non isolated zeros. The set $\bar{\Nscr}$ is at most countable, and therefore $|\bar{\Nscr} |=0$. These zeros do not contribute to the value of $\Sigma(g)$. The set $\Nscr$ is closed.
 
 \begin{definition}
 We say that $g \in H_{\br}^2$ is \emph{well behaved} if $\Nscr$ is empty or the union of finitely many disjoint closed intervals, i.e. 
 \begin{align*}
 \Nscr = \cup_{j=1}^k [\ell_j,r_j]
 \end{align*}
 for some $k \geq 1$ and $0 \leq \ell_1 < r_1 < \cdots < \ell_k < r_k \leq 1$.
 \end{definition}
 \begin{lemma}\label{le:improveFunction}
 For any $g \in H_{\br}^2$ and $\delta>0$ there exists a \emph{well behaved} function $\hat{g} \in H_{\br}^2$ such that $\| g - \hat g\|_\infty < \delta$ and $\Sigma(\hat g) \leq \Sigma(g)$.
 \end{lemma}
\begin{proofsect}{Proof}
We start by observing that for $t \in \Nscr$, we have $g'(t) =0$. Indeed, by definition there exists a sequence $(t_n)$ in $[0,1] \setminus \{t\}$ which converges to $t$ and along which $g$ vanishes. Hence
\begin{align*}
g'(t) = \lim_{n \to \infty} \frac{g(t_n) - g(t)}{t_n - t} =0.
\end{align*}
By uniform continuity of $g$ there exists a $\delta'$ such that for $|t-t'| < \delta'$ we have $|g(t) - g(t')| < \delta$. We define  recursively
\begin{align*}
\ell_1 = \inf \Nscr   \qquad r_1 &= \inf\{t \in \Nscr \colon (t,t+\delta') \cap \Nscr = \emptyset   \}, \\
\ell_2 = \inf \{ t \in \Nscr \colon t >r_1\}   \qquad r_2 &= \inf\{t \in \Nscr \colon t> \ell_2 \, \text{and } (t,t+\delta') \cap \Nscr = \emptyset   \}, 
\end{align*}
and so on. Then we set $\hat{g}=0$ on the intervals $[\ell_j,r_j]$ and $\hat{g}=g$ elsewhere. The function $\hat{g}$ constructed in this way satisfies the desired properties.
\qed
\end{proofsect} 
 Lemma \ref{le:improveFunction} implies that it suffices to establish \eqref{e:LDPlowbou} for \emph{well behaved} functions $g$ and from now on we will assume that $g$ is well behaved. 
 Furthermore, in the case where $\Nscr = \emptyset$ the bound \eqref{e:LDPlowbou} follows from the Gaussian LDP, so that we can assume $\Nscr \neq \emptyset$.
 We will first discuss the notationally simpler case where 
  $\Nscr$ consists of a \emph{single} interval $[\ell,r]$ for $0 < \ell < r < 1$. We explain how to extend the argument to the general case in the last step. 
\end{proofsect}

 \begin{proofsect}{Expansion and ``good pinning sites''}
 From now on we assume that there exist $0 < \ell < r < 1$ such that $g=0$ on $\Nscr = [\ell,r]$ and such that all zeros of $g$ outside of $\Nscr$ are isolated. 
Under these assumptions we will show that
 \begin{align}\notag
  &\liminf_{N\to\infty}\frac{1}{N} \log  \frac{\ZNe(\br)}{Z_{N}(\zero)} 
 \gNr(\|h_N -g\|_\infty < \delta)\\
   \label{e:LDPlowbou2}
 & \qquad \ge - \Big( \frac12 \int_0^\ell  \ddot{g}^2(t)\,\d t  -\tau(\eps)  (r-\ell)  + \frac12 \int_r^1  \ddot{g}^2(t)\,\d t  \Big) .
 \end{align}

The definition \eqref{gibbsdist} of $\gNr$ can be rewritten as
\begin{align}\label{gibbsdist1}
 \ZNe(\br) \;& \gNr(\d\phi)=\notag\\
 & \sum_{\Pscr \subseteq \{1, \ldots, N-1 \}}\ex^{-\Hscr_{[-1,N+1]}(\phi)} \prod_{k \in \Pscr} \eps\delta_0(\d\phi_k)   \prod_{k \in \{1, \ldots, N-1 \} \setminus \Pscr} \d\phi_k \prod_{k\in\{-1,0,N,N+1\}}\delta_{\psi_k^{\ssup{N}}}(\d\phi_k).
\end{align}
The first crucial observation is that for certain choices of ``pinning sites'' $\Pscr$ the right hand side of this expression becomes a product measure. Indeed, if $\Pscr$ contains two adjacent sites $p, p+1$ we can write
\begin{align*}
\Hscr_{[-1,N+1]}(\phi) = \Hscr_{[-1,p]}(\phi)  +\frac12 (\Delta \phi_p)^2  + \frac12 (\Delta \phi_{p+1})^2  +  \Hscr_{[p+1,N+1]}(\phi), 
\end{align*}
 which turns into $\Hscr_{[-1,p]}(\phi)  + \frac12  \phi_{p-1}^2  + \frac12  \phi_{p+2}^2  +  \Hscr_{[p+1,N+1]}(\phi)$ if $\phi_p = \phi_{p+1} = 0$. This  means that when $\phi_p$ and $\phi_{p+1}$ are pinned, the Hamiltonian decomposes into two independent contributions -- one which depends only on (the \emph{left} boundary conditions on $\phi(-1)$, $\phi(0)$  given in \eqref{bc} and)  $\phi_1, \ldots, \phi_{p-1}$ 
 and one which only depends on $\phi_{p+2}, \ldots, \phi_{N-1}$ (and the \emph{right} boundary conditions on  $\phi(N),\phi(N+1)$). Then the term corresponding to this choice of $\Pscr$ in the expansion \eqref{gibbsdist1} factorises into two independent parts. We will now reduce ourselves to choices of pinning sites $\Pscr$ which have this property.
 \begin{definition}\label{def:very-good}
 For $N \geq 2$ set $p_* := \lfloor N \ell \rfloor $ and $p^* := \lfloor N r \rfloor$.
A subset $\Pscr \subseteq \{ 1, \ldots, N-1\}$ is a \emph{very good } choice of pinning sites  if 
\begin{itemize}
\item $\{ 1, \ldots , p_*-1\} \cap \Pscr = \emptyset$ and  $\{ p^*+ 1, \ldots, N-1\} \cap \Pscr =\emptyset$.
\item $\{ p_*,\, p_* +1,\; p^* , p^* -1 \} \subseteq \Pscr$. 
\end{itemize}
(Here we leave implicit the $N$-dependence of $p_*$ and $p^*$). 
 \end{definition}
As all the terms in \eqref{gibbsdist1} are non-negative we can obtain a lower bound by reducing the sum to \emph{very good} $\Pscr$. In this way we get
 \begin{align}
  &\ZNe(\br)  \gNr(\|h_N -g\|_\infty < \delta)\notag\\
  & \geq \sum_{\Pscr \text{ very good}} \eps^{|\Pscr| } \;\ZNl(\br) \;\gNll \big(\sup_{0 \leq t \leq \ell } | h_N(t) - g(t)| \leq \delta    \big) \notag\\
  & \qquad    \times \ZNP(\zero) \,   \gNP\big(\sup_{\ell \leq t \leq r } | h_N(t) | \leq \delta    \big) \;\notag\\
  & \qquad \times Z_{[p^*-1,N+1]}(\br) \, \gamma_{[p^*-1,N+1]}^{\br} \big(\sup_{r \leq t \leq N} | h_N(t) - g(t)| \leq \delta    \big).\label{e:product_representation}
 \end{align}
The measures $\gNll $  and $ \gamma_{[p^*,N+1]}^{\br}$ on the right hand side of this expression are defined as
\begin{align*}
\gNll (\d  \phi) &=  \frac{1}{ \ZNl(\br) }\ex^{-\Hscr_{[-1,p_*+1]}(\phi)  }  \prod_{k \in \{1, \ldots, p_* -1 \} } \d\phi_k \prod_{k\in\{-1,0, p_*, p_*+1\}}\delta_{\psi_k^{\ssup{N}}}(\d\phi_k)  \\
\gamma_{[p^*-1,N+1]}^{\br}(\d  \phi) &=  \frac{1}{ Z_{[p^*-1,N+1]}(\br) }\ex^{-\Hscr_{[p^*-1,N+1]}(\phi)  } \;  \; \prod_{k \in \{p^*+1 ,\ldots, N-1 \} } \d\phi_k \prod_{k\in\{p^*-1, p^*,N,N+1\}}\delta_{\psi_k^{\ssup{N}}}(\d\phi_k) .
\end{align*}  These measures do not depend on the specific choice $\Pscr$ of very good pinning sites. The measure  $\gNP$ is defined as
\begin{align*}
\gNP(\d  \phi) =  \frac{1}{ \ZNP (\zero) }\ex^{-\Hscr_{[p_*,p^*]}(\phi) } \prod_{k \in \Pscr}  \delta_0(\d \phi_{k})  \;  \prod_{k \in \{ p_* ,  \ldots, p^* \} \setminus \Pscr } \d\phi_k .
\end{align*}
Note that none of these measures depends on the choice $\eps$ of  pinning strength, which only appears as a factor $\eps^{|\Pscr|}$ in each term in \eqref{e:product_representation}. Note furthermore, that all three measures $\gNll $, $\gamma_{[p^*-1,N+1]}^{\br} $ and $\gNP$ are Gaussian. 
\begin{lemma}\label{le:GaussLDP}
For every $\epsilon >0$ there exists an $N_*< \infty$ such that for all $N\geq N^*$  we have
\begin{align}
 \gNll\big(\sup_{0 \leq t \leq \ell } | h_N(t) - g(t)| \leq \delta    \big) &\geq \exp\Big(- N\Big[  \int_{0}^{\ell}  \ddot g(t)^2 \d t  - \inf_{h} \int_{0}^{\ell}  \ddot h(t)^2 \d t   +  \epsilon\Big] \Big) \label{e:GauLDP1} \\
 \gamma_{[p^*-1,N+1]}^{\br}\big(\sup_{r \leq t \leq 1} | h_N(t) - g(t)| \leq \delta    \big)  & \geq \exp\Big(- N \Big[ \int_{r}^{1}  \ddot g(t)^2 \d t -\inf_h \int_{r}^1  \ddot h_r(t)^2 \, \d t +\epsilon  \Big] \Big)\;, \label{e:GauLDP2}
\end{align}
where the infimum is taken over all $h\colon [0,\ell] \to \R$ and $h\colon [r,1] \to \R$ which satisfy the right boundary conditions, i.e.  $h(0) = a$, $\dot h (0) =\alpha $, $h(\ell) = 0$, $\dot h(\ell)=0$ for \eqref{e:GauLDP1} and 
$h(r) = 0$, $\dot h (r) =0 $, $h(1) = b$, $\dot h(1)=\beta$ for \eqref{e:GauLDP2}.
\end{lemma} 
 \begin{proofsect}{Proof}
This follows immediately from the Gaussian large deviation principle presented in Proposition~\ref{LDP-nopinning}.
 \qed
 \end{proofsect}
 
 \begin{lemma}\label{le:VarianceBound}
 There exists an $N_*< \infty$ such that for $N \geq N_*$ and for all very good $\Pscr \subseteq \{1, \ldots, N-1 \}$ we have
 \begin{align*}
\gNP\big(\sup_{\ell \leq t \leq r } | h_N(t) | \leq \delta    \big) \geq \frac12  \;.
 \end{align*}
 \end{lemma}
 \begin{proofsect}{Proof}
By the definition of $h_N$ we get
\begin{align*}
\gNP\big(\sup_{\ell \leq t \leq r } | h_N(t) | > \delta    \big) &\leq \gNP\big(\sup_{p_* \leq k \leq p^* } | \phi(k) | >\delta N^2    \big) \\
& \leq \sum_{p_* \leq k \leq p^* }\gNP\big( | \phi(k) | > \delta N^2    \big)\;.
\end{align*}
Recall that under $\gNP$ all $\phi(k)$ are centred Gaussian random variables and that the sum on the right hand side goes over at most $p^* - p_*+1 \leq N$ terms. Hence in order to conclude it is sufficient to prove that for all $N$ and for  all $\Pscr$ and for all $k \in \{ p_* , \ldots, p^* \}$ the variance of $\phi(k)$ under $\gNP$ is bounded by $N^3$.

To see this, we recall a convenient representation of Gaussian variances: If $C$ be the covariance matrix of a centred non-degenerate Gaussian measure on $\R^N$. Then we have for $k = 1, \ldots, N$,
\begin{align*}
C_{k,k} = \sup_{y \in \R^N \setminus \{ 0\}} \frac{y_k^2}{ \langle y, C^{-1} y \rangle} \;,
\end{align*}
where $\langle \cdot, \cdot \rangle$ denotes the canonical scalar product on $\R^N$. This identity follows immediately from the Cauchy-Schwarz inequality. 
In our context, this implies that the variance of $\phi(k)$ under $\gNP$ is given by  
\begin{align*}
\sup_{ \substack{\eta \colon \{ p_* , \ldots, p^* \} \to \R \\ \eta (k) = 0 \, \text{ for  } k \in \Pscr } } \frac{\eta(k)^2}{ 2 \, \Hscr_{[p_*,p^*]}(\eta)} \leq  \sup_{ \substack{\eta \colon \{ p_* , \ldots, p^* \} \to \R \\ \eta (k) = 0 \, \text{ for  } k \in \{p_*,p_*+1, p^*-1,p^* \} } } \frac{\eta(k)^2}{ 2 \,\Hscr_{[p_*,p^*]}(\eta)} \;,
\end{align*}
where the inequality follows because the supremum is taken over a larger set. 

The quantity on the right hand side can now be bounded easily. By homogeneity we can reduce the supremum to \emph{test vectors} $\eta$ that satisfy $\eta(k)=1$. Invoking the homogeneous boundary conditions, for such $\eta$ there must exist a $j \in \{p_* , \ldots, p^* \}$ such that $\eta(j+1) - \eta(j) \geq \frac{1}{N}$. Invoking the homogenous boundary conditions once more (this time for the difference  $\eta(p_*+1) - \eta(p_*)$) we get
\begin{align*}
\frac{1}{N}& \leq  \sum_{m = p_*+1}^{j} \big(\eta(m+1) - \eta(m) \big) - \big(\eta(m) - \eta(m-1) \big) 
= \sum_{m = p_*+1}^{j} \Delta \eta(m)  
\leq \sum_{m = p_*+1}^{p^*-1} |\Delta \eta(m) \big| \\
&\leq (p^* - p_*-1)^{\frac12} \Big( \sum_{m = p_*+1}^{p^*-1} |\Delta \eta(m) \big|^2\Big)^{\frac12}\;.
\end{align*}
Using the bound $p^* - p_*-1 \leq N$ we see that $\eta$ must satisfy  
\begin{align*}
\Hscr_{[p_*, p^*]}(\eta) \geq \frac{1}{2N^3} \;,
\end{align*}
which implies the desired bound on the variance.
\qed
\end{proofsect}
  \end{proofsect}

 \begin{proofsect}{The pinning potential}
 First of all, we observe that the minimal energy terms appearing in \eqref{e:GauLDP1} and \eqref{e:GauLDP2} can be absorbed into the boundary conditions. We obtain by  the identity \eqref{partition} in conjunction with Proposition~\ref{minimiserconv} that for every $\epsilon>0$ and for $N$ large enough 
 
  \begin{align*}
\exp\big( N  \inf_{h} \int_{0}^{\ell}  \ddot h(t)^2 \d t \big) \ZNl(\br) \geq \ZNl (\zero) \exp(-\epsilon N) \\
\exp\big( N  \inf_{h} \int_{r}^{1}  \ddot h(t)^2 \d t \big)   Z_{[p^*,N+1]}(\br) \geq Z_{[p^*,N+1]}(\zero)\exp(- \epsilon N) \;.
 \end{align*}
 Therefore, combining \eqref{e:product_representation} with Lemma~\ref{le:GaussLDP} and Lemma~\ref{le:VarianceBound} we obtain for any $\epsilon >0$ and for $N$ large enough
 $$
 \begin{aligned}
 \frac{  \ZNe(\br)}{Z_{N}(\zero)}  \gNr(\|h_N -g\|_\infty < \delta)& \geq 
  \frac12 \exp\Big(- N \int_{0}^{\ell}  \ddot g(t)^2 \d t -  N \int_{r}^{1}  \ddot g(t)^2 \d t - N\epsilon \Big) \\
  &\ \qquad \qquad   \times \sum_{\Pscr \text{ very good}} \eps^{|\Pscr|} \frac{\ZNl(\zero)   \ZNP(\zero)  \; Z_{[p^*-1,N+1]}(\zero)  }{Z_{N}(\zero)}  \;.
 \end{aligned}
$$
It remains to treat the sum of the partition functions on the right hand side. First of all, we observe that $\ZNl(\zero) $ and $ Z_{[p^*-1,N+1]}(\zero) $ and   $ Z_N(\zero) $ do not depend on the choice of very good  $\Pscr$ so that they can be taken out of the sum, i.e. we can write 
\begin{align*}
&\sum_{\Pscr \text{ very good}} \eps^{|\Pscr|} \frac{\ZNl(\zero)   \ZNP(\zero)  \;\ZNrr(\zero)  }{Z_{N}(\zero)}  \\
&\qquad =  \frac{\ZNl(\zero)  Z_{[p^* , p_*]}(\zero)  \;\ZNrr(\zero)  }{Z_{N}(\zero)} 
 \sum_{\Pscr \text{ very good}} \eps^{|\Pscr|} \frac{ \ZNP(\zero)  }{Z_{[p^* , p_*]}(\zero)} \;.
\end{align*}
Here we have multiplied and divided by the Gaussian partition function $Z_{[p^* , p_*]}(\zero)$ (In the notation of the introduction this constant could also be written as  $Z_{p^* -p_* -2}$, but we prefer to keep the explicit dependence on the interval  in the notation). This  allows us to compare the sum on the right hand side to the limit  \eqref{freeenergies} which defines $\tau(\eps)$. More precisely we get
\begin{align*}
 \sum_{\Pscr \text{ very good}} \eps^{|\Pscr|} \frac{ \ZNP(\zero)  }{Z_{[p^*, p_*]}(\zero)}  = \frac{Z_{[p^*, p_*],\eps}(\zero)}{Z_{[p^*, p_*]}(\zero)} \geq \exp\big((r-\ell) \tau(\eps)  - N \epsilon \big)  \;
\end{align*}
for $N$ large enough (depending on $\epsilon$), where the equality follows from reversing the expansion. 
To conclude it only remains to observe that according to Appendix~\ref{AppB} the quotient
\begin{align*}
 \frac{\ZNl(\zero)  Z_{[p^*, p_*]}(\zero)  \;\ZNrr(\zero)  }{Z_{N}(\zero)} 
 \end{align*}
 decays at most polynomially in $N$ which implies that it disappears on an exponential scale. Therefore, \eqref{e:LDPlowbou2} follows. 

\end{proofsect}
We have thus established \eqref{e:LDPlowbou} for an open ball around a well behaved function which has exactly one zero interval. As outlined earlier after Lemma~\ref{le:improveFunction} we actually need to show  \eqref{e:LDPlowbou} for all well behaved functions. For a general well behaved functions with $ \Nscr = \cup_{j=1}^k [\ell_j,r_j]$ and $0 \leq \ell_1 < r_1 < \cdots < \ell_k < r_k \leq 1$ the proof can be easily adapted: For $j = 1, \ldots, k$ we define the discrete boundary points $p_{*,j} = \floor{N\ell_j}$ and $p_j^{*}= \floor{N r_j}$ and define very good pinning sites   to be those subsets of $\{1,\ldots, N-1\}$  which  contain all of the $p_{*,j}, p_{*,j}+1, p^{*}_j-1, p^{*}_j$  and none of the sites to the left of $p_{*,1}$, between the  $p_{j}^{*} $ and $p_{*,j+1}$,  or to the right of $p^{*}_k$. In the product representation \eqref{e:product_representation} we then get a larger number of independent factors -- one for each of the $k$ pinned intervals and one for each of the $k+1$ intervals where the interface can move away from the $x$-axis (in the case where $\ell_1 = 0$ or $r_k = 1$ there are only $k$ or even $k-1$  intervals where the interface can move away). Lemma~\eqref{le:GaussLDP} can then be applied to each of the ``free'' intervals and Lemma~\ref{le:VarianceBound} can be applied to each of the ``pinned'' intervals and the discussion of the partition functions can be repeated with only obvious changes.

Finally we mention that  the case of Dirichlet boundary conditions on the left hand side and free boundary conditions on the right hand side follows in the exact same way. The only difference is that the right boundary condition in the definition of $\gamma_{[p^*,N+1]}^{\br}$ should be removed and that consequently the infimum in \eqref{e:GauLDP2} has to be taken over the larger class of all $h$ satisfying $h(r) = \dot{h}(r) =0$ without any restriction on $h(1)$ or $\dot{h}(1)$.

\subsubsection{Proof of the upper  bound of Theorem~\ref{THM-main}}\label{mainupper}  
For the upper bound we need to show that 
\begin{equation}\label{e:up1}
\limsup_{N \to \infty} \frac{1}{N} \log \frac{\ZNe(\br)}{Z_N(\zero)}  \gNr (h_N \in\Kscr) \leq - \inf_{h \in \mathcal{\Kscr}} \Sigma(h),
\end{equation}
for all closed $\Kscr \subset \Cscr([0,1];\R) $. 
 \begin{proofsect}{Reduction to a simpler statement}
First of all we observe:
\begin{lemma}\label{le:tightness}
For any $N\in\N$ let $ \gNr$ be the measure  given in \eqref{gibbsdist} with boundary conditions as in \eqref{bc} and let the rescaled profiles $h_N$ be as given in \eqref{def:hN}.  Then the sequence of 
 distributions of the rescaled profiles  $h_N$ is exponentially tight in $ \Cscr([0,1];\R) $.
\end{lemma}
The proof of this lemma can be found at the end of this section. Lemma~\ref{le:tightness} implies that it suffices to establish \eqref{e:up1} for compact sets $\Kscr$. %
Going further, it suffices to show that for any $g \in \Cscr([0,1];\R) $ and any $\epsilon>0$ there exists a $\delta =\delta(g,\epsilon)>0$ such that 
\begin{align}\label{e:up3}
\limsup_{N \to \infty} \frac{1}{N} \log \frac{\ZNe(\br)}{Z_N(\zero)} \gNr(h_N \in B(g, \delta) ) \leq -  \Sigma (g)+\epsilon .
\end{align}
Here $B(g,r) = \{ h \in \Cscr([0,1];\R)  \colon \| h - g\|_\infty < r\}$  denotes the $L^\infty$ ball of radius $r$ around $g$. 

We give the 
simple argument to show that \eqref{e:up3} implies \eqref{e:up1}: For any compact set $\Kscr$ and any $\epsilon>0$ there exists a finite set $\{ g_1, \ldots, g_M\} \subset \Kscr$ such that $\Kscr \subseteq \cup_{j=1}^M B(g_j, \delta(g_j, \epsilon))$. Then \eqref{e:up3} yields
\begin{align*}
\limsup_{N \to \infty} \frac{1}{N} \log \frac{\ZNe(\br)}{Z_N(\zero)} \gNr (h_N \in \Kscr)& \leq  \limsup_{N \to \infty} \frac{1}{N} \log \sum_{j=1}^M \frac{\ZNe(\br)}{Z_N(\zero)} \gNr (h_N \in B(g_j, \delta(g_j,\epsilon)) )\\
& \leq \max_{j=1, \ldots, M} \limsup_{N \to \infty} \frac{1}{N} \log  \frac{\ZNe(\br)}{Z_{N}(\zero)} \gNr (h_N \in B(g_j, \delta(g_j,\epsilon) ) )\\
& \leq -\min_{j=1, \ldots, M}  \Sigma (g_j)+ \epsilon \\
& \leq - \inf_{h \in \Kscr}  \Sigma (h) +\epsilon ,
\end{align*}
so \eqref{e:up1} follows because $\epsilon>0$ can be chosen arbitrarily small.

For fixed $g$ and $\epsilon$ the value of $\delta$ is determined by the following lemma.
\begin{lemma}\label{lem:very-well-behaved}
For any  $g \in\Cscr([0,1];\R)$  and  all $\epsilon>0$ there exists a $\bar{\delta}>0$ and a closed set $\mathcal{I} \subset [0,1]$
such that the following hold:
\begin{enumerate}
\item $\mathcal{I}$ is the union of finitely many disjoint closed intervals, i.e. 
\begin{equation}
\mathcal{I} = \cup_{j=1}^M [\ell_j, r_j]
\end{equation}
for some finite $M$ and $0 \leq  \ell_1 < r_1< \ell_2 < r_2 < \ldots < r_M \leq 1$.
\item The level-set $\{ t \in [0,1] \colon |g(t) | \leq \bar{\delta} \}$ is contained in $\mathcal{I}$.
\item The measure of $\mathcal{I}$ satisfies the bound 
\begin{equation*}
| \mathcal{I}| \leq | \{ t \in [0,1] \colon |g(t) | = 0 \}|+\epsilon.
\end{equation*}
\end{enumerate}
\end{lemma}
The proof of this lemma is also given at the end of this section.
\end{proofsect}

\begin{proofsect}{Expansion and key lemmas}
We will now proceed to prove that \eqref{e:up3} holds for a fixed $g$ and $\epsilon$ and a suitable $\delta \in (0, \bar{\delta})$, where $\bar{\delta}$ is given by Lemma~\ref{lem:very-well-behaved}.
For simplicity (and  similar to the proof of the lower bound), we will assume that the set $\Ical$ constructed in this lemma consists of \emph{a single} interval $[\ell,r]$. The argument for the  case of a finite union of disjoint intervals  is identical, only requiring slightly more complex notation, and will be omitted.

\medskip
We write
\begin{align}\label{e:up4}
\gNr (h_N \in B(g, \delta)) = \frac{1}{\ZNe(\br) } \sum_{\Pscr \subseteq \{ 1, \ldots, N-1\} } \eps^{|\Pscr|} \;\ZNPP(\br) \; \gNPP(h_N \in B(g,\delta)),
\end{align}
where as above $\gNPP$ denotes the Gaussian measure over $\{-1, \ldots, N+1 \}$ which is pinned at the sites in $\Pscr$, i.e. 
\begin{align*}
\gNPP (\d  \phi) &=  \frac{1}{ \ZNPP(\br) }\ex^{-\Hscr_{[-1,p_*+1]}(\phi)  }  \prod_{k \in \{1, \ldots, N-1 \} \setminus \Pscr } \d\phi_k \prod_{k\in \Pscr \cup \{-1,0, N, N+1\}}\delta_{\psi_k^{\ssup{N}}}(\d\phi_k) , 
\end{align*}
and $\ZNPP(\br)$ is the corresponding Gaussian normalisation constant.
By definition $|g(t) | >\bar{\delta}>\delta$ for $t \in [0,1] \setminus \Ical$, so in \eqref{e:up4} it suffices to sum over those sets of pinning sites $\Pscr \subseteq N \Ical \cap \Z$. The next two lemmas simplify the expressions in the sum \eqref{e:up4}. For the moment we only deal with homogeneous boundary conditions $\br=0$ and start by introducing some notation which will be used to simplify the  partition functions.
As above in Definition~\ref{def:very-good} we will be interested in sets of pinning sites $\Pscr$ that allow to separate the Hamiltonian $\Hscr_{[-1,N+1]}$ into independent parts.
\begin{definition}
Let $\Pscr \subseteq \{ 1 ,\ldots, N-1 \}$ be non-empty and let $p_* = \min \Pscr$ and $p^*  = \max \Pscr$. We will call $\Pscr$ an \emph{good} choice of pinning sites if 
$\{ p_*+1, p^*-1 \} \subseteq \Pscr$. 
We will also call the empty set good. 
\end{definition}
Note that the \emph{very good} sets introduced in Definition~\ref{def:very-good} are  \emph{good} but the inverse implication is is not true. The difference between the 
two notions is that we do not prescribe the precise value of $p_*$ and $p^*$ for good sets. They will however always be confined to the 
interval $[ \lfloor \ell N \rfloor, \lfloor r N \rfloor +1]$. 
We also introduce the following operation of \emph{correcting} a set to make it good. 
\begin{definition}
Let $\Pscr \subseteq \{ 1, \ldots, N-1 \}$ be non-empty with $p_* = \min \Pscr$ and $p^*  = \max \Pscr$. Then we define 
\begin{equation*}
c(\Pscr) = \Pscr \cup \{ p_*+1, p^*-1\}.
\end{equation*}
We also set $c(\emptyset) = \emptyset$.
\end{definition}
For later use we remark that on the one hand the correction map $c$ adds at most two points to a given set $\Pscr$, and that on the other hand for a given good set $\Pscr$ there are at most $4$ distinct $\widetilde{\Pscr}$ with $c(\widetilde{\Pscr}) = \Pscr$.
The following Lemma permits to replace the partition function $\ZNPP(\zero) $ in \eqref{e:up4} by the partition function  $Z_{[-1,N+1] \setminus c(\Pscr)}(\zero)$ with corrected choice of pinning sites. 
\begin{lemma}\label{le:correction1}
For every non-empty $\Pscr \subseteq \{ 1 , \ldots, N-1 \}$ we have 
\begin{align*}
\ZNPP(\zero) \leq (2 \pi N) Z_{[-1,N+1] \setminus c(\Pscr)}(\zero).
\end{align*}
\end{lemma}
\begin{proofsect}{Proof}
For any $\Pscr \subset \{ 1, \ldots, N-1 \}$ set
\begin{align*}
\gamma_{[-1,N+1] \setminus \Pscr}^{\zero}(\d  \phi) &=  \frac{1}{ Z_{[-1,N+1] \setminus \Pscr}(\zero) }\ex^{-\Hscr_{[-1,N+1]}(\phi)  }  \prod_{k \in \{1, \ldots, N-1 \} \setminus \Pscr } \d\phi_k 
 \prod_{k\in  \{-1,0,N, N+1\} \cup \Pscr}\delta_{0}(\d\phi_k) .
\end{align*}
We derive an identity that links the Gaussian partition function $Z_{[-1,N+1] \setminus \Pscr}(\zero)$ to $Z_{[-1,N+1] \setminus (\Pscr \cup \{ j\})  } (\zero)$ for an arbitrary $\Pscr \subseteq \{1, \ldots, N-1\}$ and $j \notin \Pscr$. We have
\begin{align}
\notag
Z_{[-1,N+1] \setminus \Pscr}(\zero)  = &\int_\R \Big( \int  \ex^{-\Hscr_{[-1,N+1]}(\phi)} \prod_{k \in \{1, \ldots, N-1 \} \setminus (\Pscr \cup \{ j\}) } \d\phi_k \\
\label{e:up-lemma1}
& \qquad \times \prod_{k\in  \{-1,0,N, N+1\} \cup \Pscr} \delta_{0}(\d\phi_k) \Big)  \;\d \phi_j.
\end{align}
Denote by $ \phi^*$  the unique minimiser of $\Hscr_{[-1,N+1]}$ subject to the constraints that $\phi^*(k)=0$ 
for $k  \in \{ -1 ,0,N,N+1|\} \cup \Pscr$ and $\phi^*(j)= 1$. Then by homogeneity for any $y \in \R$ the function $y \phi^*(k)$
is the unique minimiser of $\Hscr_{[-1,N+1]}$ constrained to be zero on the same set, but satisfying $y \phi^*(j) =y$.
This implies that for any  $\phi \colon \{ -1, \ldots, N+1 \} \to \R$ satisfying the same pinning constraint we have
\begin{align*}
\Hscr_{[-1,N+1]}(\phi) = \Hscr_{[-1,N+1]}(\phi- \phi(j) \phi^*) + \phi(j)^2 \Hscr_{[-1,N+1]}( \phi^*). 
\end{align*}
As in the proof of Lemma~\ref{le:VarianceBound} we can see that $\Hscr_{[-1,N+1]}( \phi^*) = \frac{1}{2\mathrm{var} (\phi(j))}$ where 
$\mathrm{var}(\phi(j))$ denotes the variance of $\phi(j)$ under $\gamma_{[-1,N+1] \setminus \Pscr}^{\zero}$. 
This allows to rewrite \eqref{e:up-lemma1} as 
\begin{align}
Z_{[-1,N+1] \setminus \Pscr}(\zero)  =& \int_\R \Big( \int \ex^{-\Hscr_{[-1,N+1]}(\phi- \phi(j) \phi^*)} \prod_{k \in \{1, \ldots, N-1 \} \setminus (\Pscr \cup \{ j\}) } \d\phi_k \\
& \qquad \times \prod_{k\in  \{-1,0,N, N+1\} \cup \Pscr}\delta_{0}(\d\phi_k)  \Big) \ex^{-\frac{y^2}{2\mathrm{var} (\phi(j)) }} \d y\notag\\
=& Z_{[-1,N+1] \setminus (\Pscr \cup \{ j\})  } (\zero)\int_{\R} \ex^{\frac{-y^2}{2\mathrm{var} (\phi(j)) }} \d y = Z_{[-1,N+1] \setminus (\Pscr \cup \{ j\})  } (\zero) \sqrt{2 \pi \,\mathrm{var} (\phi(j))}\notag.
\end{align}
As the correction map $c$ adds at most two points to the pinned set it only remains to get an upper bound on the variance of $\phi(j)$ for $j = p_* +1$ or $j = p^*-1$ under $\gamma_{[-1,N+1]\setminus \Pscr}^{\zero}$, or equivalently a lower bound on $\Hscr_{[-1,N+1]}( \phi^*)$; we show the argument for $p_*+1$. It is very similar to the upper bound on the variance derived in Lemma~\ref{le:VarianceBound}, but this time we obtain a better bound using the fact that $p_*+1$ is adjacent to a pinned site. More precisely, using that $\phi^*(p_*)=0$ and  
the fact that the homogenous boundary conditions imply $\phi^*(0) - \phi^*(-1) =0$, we get
\begin{align*}
1& = \phi^*(p_*+1) - \phi^*(p_*) = \sum_{j=0}^{p_*} \big( \phi^*(j+1) - \phi^*(j) \big) - \big( \phi^*(j) - \phi^*(j-1) \big) = \sum_{j=0}^{p_*}  \Delta \phi^*(j)\\
&\leq (p_*+1)^{\frac12}  \Big(\sum_{j=0}^{p_*}  (\Delta \phi^*(j))^2 \Big)^{\frac12} \leq N^{\frac12}\big( 2 \Hscr_{[-1,N+1]}( \phi^*) \big)^{\frac12}.
\end{align*}
This finishes the argument.
\qed
\end{proofsect}
The next Lemma provides an upper bound on the Gaussian probabilities appearing in \eqref{e:up4} (still for homogeneous boundary conditions). It is 
essentially a variant of the Gaussian isoperimetric inequality.  To state it, we introduce the rescaled Hamiltonian
\begin{equation*}
\Escr_N(h) = \frac12 \sum_{j=0}^N  \frac{1}{N}    N^4  \Big( h\Big(\frac{j+1}{N} \Big) + h\Big(\frac{j-1}{N} \Big) - 2 h\Big(\frac{j}{N} \Big) \Big)^2.
\end{equation*}
Observe that for $h$ and $\phi$ related by \eqref{def:hN} we have
\begin{align*}
 \Hscr_{[-1,N+1]}(\phi) = N  \Escr_N(h_N).
\end{align*}

\begin{lemma}\label{le:Isoperi}
 For every $\delta >0$ there exists  an $N_0>0$ such that for all $N \geq N_0$ and all $\Pscr \subseteq \{ 1, \ldots, N-1 \}$ and all  $g \in  \Cscr([0,1];\R) $.
\begin{align*}
\gamma_{[-1,N+1] \setminus \Pscr}^{\zero} (h_N \in B(g,\delta)) \leq \exp \big(-N \inf_h \Escr_N(h)  \big), 
\end{align*}
where the infimum is taken over all $h \colon \{-\frac{1}{N}, 1, \ldots, 1, 1 + \frac{1}{N} \} \to \R$ with $\| h - g \|_{\infty} \leq 2 \delta$.
\end{lemma}
\begin{proofsect}{Proof}
We recall a convenient version of the Gaussian isoperimetric inequality (see e.g. \cite{Led}): Let $\gamma$ be a centred Gaussian measure on $\R^N$ with 
Cameron-Martin norm $\| \cdot \|_{CM}$. Furthermore, let $B \subseteq \R^N$  be a closed set satisfying $\gamma (B) \geq \frac12$. Then 
\begin{align}\label{e:IsoPP}
\gamma\big( B + \Sscr(r)  \big) \geq  \Phi(r).
\end{align}
Here 
\begin{align*}
\Phi( r) &= \frac{1}{\sqrt{2 \pi}} \int_{-\infty}^r \ex^{-\frac{x^2}{2}}\, \d x  
\end{align*}
denotes the  distribution function of the standard normal distribution, $\Sscr(r) = \{ x \in \R^N \colon \| x \|_{CM} \leq r \}$ is a closed ball in the Cameron-Martin norm 
$$ B + \Sscr(r) := \{ x + y \colon x \in B \text{ and } y \in \Sscr(r) \}.
$$

We apply this theorem to the distribution of the rescaled profile $h_N$ under the measures $\gamma_{[-1,N+1] \setminus \Pscr}^{\zero}$ . All of these distributions are Gaussian and for each choice of $\Pscr$ the Cameron-Martin norm is given by 
$\sqrt{2 N\Escr_N(h)}$ restricted to  $\R^{\{1, \ldots, N\} \setminus \Pscr}$. First of all, we can see as in Lemma~\ref{le:VarianceBound} that  for any $\delta>0$ and there exists  $N_0$ such that for $N > N_0$ we have uniformly over the choice of $\Pscr$
 \begin{align}\label{aim1}
 \gamma_{[-1,N+1] \setminus \Pscr}^{\zero}(h_N\in B(0,\delta))\geq \frac12.
 \end{align}
Indeed, just like in the proof of this Lemma, the probability of the complement goes to zero  for large $N$ uniformly over $\Pscr$, because the variances of each $h_N(k/N)$ for $k = 1, \ldots, N-1$ are bounded by $N^{-1}$ independently of the choice of $\Pscr$.
Now we invoke  \eqref{e:IsoPP} for $B = B(0, \delta)$  and observe that the ball $B(g, \delta)$
is contained in the complement of $B+\Sscr(r)$ if 
\begin{align}
r = \inf_{h \in B(g, 2 \delta)} \sqrt{2N \Escr_N(h)} .
\end{align}
This yields 
\begin{align*}
\gamma_{[-1,N+1] \setminus \Pscr}^{\zero}(h_N \in B(g, \delta) ) \leq (1 - \Phi(r)) \leq  \ex^{-\frac{r^2}{2}}.
\end{align*}
The claim then  follows from rewriting
\begin{align*}
\frac{\big( \inf_{h \in B(g, 2 \delta)} \sqrt{2N \Escr_N(h)} \big)^2}{2} = N \inf_{h \in B(g, 2 \delta)} \Escr_N(h).
\end{align*}

%
%
%
%
 \qed
 \end{proofsect}
\end{proofsect}
\begin{proofsect}{Conclusion}
We now apply these two Lemmas to the terms appearing in the sum \eqref{e:up4}. For each $\Pscr \neq \emptyset$ we can write
\begin{align}
\notag
\ZNPP(\br)& \; \gNPP(h_N \in  B(g,\delta)) \\
&= \ex^{-\Hscr_{[-1,N+1]}(\phi_{\br,\Pscr}^*)}\ZNPP(\zero) \; \gamma_{[-1,N+1] \setminus \Pscr}^{\zero} (B(g- h_{\br,\Pscr}^{*},\delta)), \label{NEW1}
\end{align}
where we have used \eqref{partition} to include the boundary conditions into the Gaussian partition function. The function $\phi_{\br,\Pscr}^*$ is the minimiser of $\Hscr_{\L_N}$ subject to the boundary conditions $\br$ and pinned 
on the sites in $\Pscr$. The profile $h_{\br,\Pscr}^{*}$ is the rescaled version of $\phi_{\br,\Pscr}^*$ and in particular $\Hscr_{[-1,N+1]}(\phi_{\br,\Pscr}^*) = N \Escr_N (h_{\br,\Pscr}^{*})$. First of all, Lemma~\ref{le:Isoperi} allows to bound for $N$ 
large enough uniformly over $\Pscr$
\begin{align*}
\gamma_{[-1,N+1] \setminus \Pscr}^{\zero} (B(g- h_{\br,\Pscr}^{*},\delta)) &\leq \exp(- N\inf_{h \in B(g- h_{\br,\Pscr}^{*},2\delta)} \Escr_N(h))\\
&= \exp\big(- N\inf_{h \in B(g,2\delta)} \Escr_N(h) + N \Escr_N( h_{\br,\Pscr}^{*}) \big),
\end{align*}
and the last term in the exponent exactly cancels the first term on the right hand side of  \eqref{NEW1}. Plugging this into the left hand side  and then using Lemma~\ref{le:correction1} yields
\begin{align*}
\ZNPP(\br)& \; \gNPP(h_N \in  B(g,\delta)) \\
&\leq \ZNPP(\zero) \;  \exp\big(- N\inf_{h \in B(g,2\delta)} \Escr_N(h) \big) \\
& \leq (2 \pi N) Z_{[-1,N+1] \setminus c(\Pscr)}(\zero) \exp\big(- N\inf_{h \in B(g,2\delta)} \Escr_N(h) \big).
\end{align*}
Finally, we claim that there exists a $\delta \in (0, \bar{\delta})$ such that for $N$ large enough
\begin{align*}
\inf_{h \in B(g,2\delta)} \Escr_N(h)  \geq \Escr(g) - \epsilon.
\end{align*}
Indeed, if this is not the case, then there exists a sequence $\delta_n \to 0$ and a sequence $N(n) \to \infty$ such that 
\begin{align*}
\| h_n  - g \|_{\infty} \leq \delta_n \qquad \text{and} \qquad  \Escr_N(h_n) < \Escr(h)-\epsilon 
\end{align*}
which contradicts Lemma~\ref{Gammaconv}. This finally allows to write for this $\delta$ and for $N$ large enough 
\begin{align*}
\ZNPP(\br)& \; \gNPP(h_N \in  B(g,\delta)) \leq  (2 \pi N) Z_{[-1,N+1] \setminus c(\Pscr)}(\zero) \exp\big(- N (\Escr(g) -\epsilon) \big).
\end{align*}
Plugging this into \eqref{e:up4} we obtain
\begin{align*}
\frac{\ZNe(\br)}{Z_{N}(\zero)} \gNr (B(g, \delta) ) \leq (2 \pi N) \exp \big(-N(\Escr(g) - \epsilon ) \big)  \sum_{\Pscr \subseteq\{  \floor{\ell N}, \ldots,  \floor{r N}+1 \}} \eps^{|\Pscr|} \frac{ Z_{[-1,N+1] \setminus c(\Pscr)}(\zero)}{Z_{N}(\zero)} .
\end{align*}
For $\Pscr = \emptyset$ we have $\frac{ Z_{[-1,N+1]\setminus c(\Pscr)}(\zero) }{Z_{N}(\zero)}=1$ by definition so this term is of lower order.
The sum over all non-empty $\Pscr$  can then  be rewritten as
\begin{align*}
&\sum_{\substack{ \Pscr \subseteq\{  \floor{\ell N}, \ldots,  \floor{r N}+1 \} \\\Pscr \neq \emptyset}} \eps^{|\Pscr|} \frac{ Z_{[-1,N+1]\setminus c(\Pscr)}(\zero) }{Z_{N}(\zero)}\\
&  \leq 4 \sum_{\substack{\Pscr \text{ good}\\ \ell N \leq p_* < p^* \leq r N }} \eps^{|\Pscr|}  \frac{ \ZNPP(\zero)}{Z_N(\zero)} 
\leq 4 \sum_{\ell N \leq k_1 < k_2 \leq r N}  \sum_{\substack{\Pscr  \text{ good} \\ p_{*} = k_1 \\ p^{*}= k_2} } \eps^{|\Pscr|}  \frac{ \ZNPP(\zero)}{Z_{N}(\zero)}\\
& = 4  \sum_{\ell N \leq k_1 < k_2 \leq r N}  \sum_{\substack{\Pscr \text{ good} \\ p_{*} = k_1 \\ p^{*}= k_2} }  \eps^{|\Pscr|}\frac{\ZNl(\zero) Z_{[p_*, p^*]\setminus \Pscr}(\zero)  \;\ZNr(\zero)  }{Z_{N}(\zero)} \\
  & = 4  \sum_{\ell N \leq k_1 < k_2 \leq r N} \frac{Z_{[-1, k_1+1]}  (\zero)  Z_{[k_1, k_2]}(\zero)  \;Z_{[k_2,N+1]}(\zero)  }{Z_{N}(\zero)}   \frac{Z_{[k_1 , k_2],\eps}(\zero)}{Z_{[k_1 ,k_2]}(\zero)}   .
 \end{align*}
We now  bound this expression by
\begin{align*}
\leq 4 N^2 \sup_{\ell N \leq k_1 \leq k_2 \leq r N} \frac{Z_{[-1, k_1+1]}  (\zero)  Z_{[k_1, k_2]}(\zero)  \;Z_{[k_2,N+1]}(\zero)  }{Z_{N}(\zero)}  \sup_{\ell N \leq k_1 \leq k_2 \leq r N} \frac{Z_{[k_1 , k_2],\eps}(\zero)}{Z_{[k_1 ,k_2]}(\zero)}
\end{align*}
and observe that according to Appendix~\ref{AppB} the first supremum grows at most polynomially in $N$ while Proposition~\ref{existenceOfLimit} implies that for $N$ large enough the second supremum is bounded by $\leq \exp (N(r-\ell) (\tau(\eps) + \epsilon))$
which establishes \eqref{e:up3}.

\end{proofsect}
\begin{proofsect}{Proofs of Lemmas}
\begin{proofsect}{Proof of Lemma~\ref{le:tightness}}
Due to the Arzel\`a-Ascoli Theorem and the fixed boundary conditions it suffices to show that
\begin{align}\label{e:tight1}
\limsup_{M \to \infty }\limsup_{N \to \infty} \frac{1}{N} \log \big( \gNr \big( h \in \Cscr([0,1];\R)  \colon \| \dot{h} \|_{\infty} \geq M \big) \big) = -\infty.
\end{align}
We recall, that according to \eqref{gibbsdist} the measure $\gNr$ can be represented as a convex combination of Gaussian measures via
\begin{align*}
\gamma_{N,\eps}^{\br}&=\sum_{\Pscr \subseteq \{1, \ldots, N-1 \}}
\eps^{|\Pscr|}  \frac{ Z_{[-1,N+1] \setminus \Pscr} (\br)}{\ZNe(\br)} \gamma_{[-1, N+1] \setminus \Pscr}^{\br},
\end{align*}
where as before $ \gamma_{[-1, N+1] \setminus \Pscr}^{\br}$ is the Gaussian measure which is determined by the energy functional $\Hscr_{[-1,N+1]}$ and the boundary conditions $\br$ as 
well as the pinning sites $\Pscr$, and $Z_{[-1,N+1] \setminus \Pscr} (\br)$ is the corresponding partition function. To show that \eqref{e:tight1} holds,  we introduce a notion of $M$-typical 
sets of pinning sites for every $M$ below. Roughly speaking  $\Pscr \subset \{1, \ldots, N-1 \}$ is $M$-typical if it does not contain any point whose distance to the boundary is $\sim M^{-\frac13} N$. The bound
\eqref{e:tight1} then follows from the following two statements:
\begin{itemize} 
\item For every choice of boundary conditions $\br =(a,\alpha,b,\beta)$ and any $M \geq 1$ there exists an $N_0$  such that for $N \geq N_0$ and for any $M$-typical $\Pscr \subset \{ 1, \ldots, N-1 \}$ we have
\begin{align}\label{e:Tight2}
 \gamma_{[-1, N+1] \setminus \Pscr}^{\br}   \big( h \in \Cscr([0,1];\R)  \colon \| \dot{h} \|_{\infty} \geq M \big)  \leq  \ex^{-\frac{NM^2}{8}}  ,
\end{align}
\item For choice of boundary conditions $\br =(a, \alpha,b,\beta )$ there exists a $c>0$  such that for any $M \geq 1$ there exists $N_0>0$  such that for $N \geq N_0$
\begin{align}\label{e:Tight3}
\sum_{\Pscr \text{ not $M$-typcial}} \eps^{|\Pscr|}  \frac{ Z_{[-1,N+1] \setminus \Pscr} (\br)}{\ZNe(\br)}  \leq   \ex^{-cN M^{\frac13}}  .
\end{align}
\end{itemize}

The rest  of the proof is devoted to  establishing the bounds \eqref{e:Tight2} and \eqref{e:Tight3}. We start with the following two Lemmas which summarise useful properties of 
the Hamiltonian $\Hscr_{[-1,N+1]}$ on configurations pinned close to the boundary. The first Lemma gives a lower bound on $\Hscr_{[-1,N+1]}(\phi)$ for profiles $\phi$ pinned close to the 
boundary. The second Lemma asserts the existence of a profile $\tilde\phi^*$  which satisfies the boundary conditions $\br$ and the pinning condition at 
the sites in $\Pscr \subset \{ 1, \ldots, N-1\}$ with a good control on $\Hscr_{[-1,N+1]}(\tilde\phi^*)$,  provided that  $\Pscr$ does not contain sites close to the boundary. The proofs of both Lemmas are given below, but before we conclude the 
proof  Lemma~\ref{le:tightness} assuming that they hold.
\begin{lemma}\label{le:Energytypical1}
Let $(a,  \alpha) \neq (0,0)$. Then there exists $\delta_0 >0$ and a  $c >0$ such that for all $N $, all $L \leq \delta_0 N$,  all $\Pscr \subset \{ 1, \ldots, N-1 \}$ with $\min \Pscr < L$ and all 
$\phi \colon \{-1, \ldots, N+1 \} \to \R$ satisfying the boundary conditions  $\phi(-1) = N^2 a- N \alpha $, $\phi(0) = N^2 a  $ as well as the pinning condition $\phi(k ) =0$ for $k \in \Pscr$ we have
\begin{align*}
\Hscr_{[-1, N+1]}(\phi) \geq c\frac{N^2}{L}.
\end{align*}
\end{lemma}
\begin{lemma}\label{le:Energytypical2}
Let $(a, \alpha) \in \R^2$. Then there exist a constant $ c$ such that for any $0<\delta<\frac12$ there exists $N_0$ such that for $N \geq N_0$ 
there exists a function $\tilde \phi^* \colon \{-1, \ldots, N+1 \} \to \R$  which satisfies the boundary conditions \eqref{bc} as well as $\tilde \phi^* (p) =0$ for 
all $\delta N \leq p \leq N - \delta N$ 
and such that 
\begin{align*}
\Hscr_{[-1, N+1]}(\tilde \phi^*) \leq c \frac{N}{\delta^3}. 
\end{align*}
\end{lemma}

Motivated by these two Lemmas we now present the definition of \emph{typical} choice of pinning sites $\Pscr$. 
\begin{definition}
Let $\br =(a, \alpha, b,\beta) \in \R^4$ and let $M \geq1$ and $N\geq 2$. Furthermore, let $\delta_0>0$ be a constant whose precise value depends on $\br$ and will be given below.
For $(a, \alpha) \neq (0,0)$ a subset $ \Pscr \subseteq \{1, \ldots, N-1 \}$  is called $M$-typical from the left if $\Pscr \cap [0, (\delta_0 M)^{-\frac13} N] = \emptyset$. For $(a, \alpha) = (0,0)$ 
and set is $M$-regular from the left.
Similarly it is called $M$-typical from the right 
if $\Pscr \cap [N(1-  (\delta_0 M)^{-\frac13} ), N] = \emptyset$ for $(b,\beta) \neq (0,0$ and any set is $M$-regular from the right for $(b,\beta)=(0,0)$. The set is $M$-typical if it is both $M$-typical from the left and from the right.
\end{definition}

We now proceed to deriving the bound \eqref{e:Tight2} for typical $\Pscr$. 
For any fixed choice of pinning sites $\Pscr \subseteq \{1, \ldots, N-1\}$  we have
\begin{align*}
& \gamma_{[-1, N+1] \setminus \Pscr}^{\br}  \big( h \in  \Cscr([0,1];\R)   \colon \| \dot{h} \|_{\infty} \geq M \big) \\
 & =  \gamma_{[-1, N+1] \setminus \Pscr}^{\br}  \big( \phi  \colon  |\phi(k+1) - \phi(k)|  \geq  N M \text{ for at least one } k \in \{ 0 , \ldots, N-1 \} \big) \\
& \leq  \sum_{k=0}^{N-1}  \gamma_{[-1, N+1] \setminus \Pscr}^{\br}  \big( \phi  \colon  |\phi(k+1) - \phi(k)|  \geq  N M \big) \\
& \leq  N \sup_{k\in \{0, \ldots, N-1 \}}  \exp \Big( -\frac{(NM- m(k))^2}{2\sigma(k)^2} \Big)    
\end{align*}
where $m(k)$ and $\sigma(k)^2$ are the mean and variance of $\phi(k+1) - \phi(k)$ under $ \gamma_{[-1, N+1] \setminus \Pscr}^{\br} $. We can thus conclude if 
we can establish that for any $M$ large enough and every $N$ large enough (depending on $M$) uniformly over all $M$-typical  $\Pscr$ and for all $k \in \{ 0, \ldots, N-1 \}$ we have
\begin{align*}
m(k) \leq  \frac12 NM \qquad \text{and}  \qquad \sigma(k)^2 \leq  N .
\end{align*}
The second bound follows by a similar argument as Lemma~\ref{le:VarianceBound} which does not make use of any specific requirements on $\Pscr$. As in this Lemma we  see
that 
\begin{align*}
\sigma(k)^2 &= \sup_{ \substack{\eta \colon \{ -1, \ldots, N+1 \} \to \R \\ \eta (k) = 0 \, \text{ for  } k \in \Pscr \cup \{ -1,0, N,N+1\} } } \frac{(\eta(k+1)- \eta(k))^2}{ 2 \, \Hscr_{[-1,N+1]}(\eta)}\\
& \leq   \sup_{ \substack{\eta \colon \{ -1, \ldots, N+1 \} \to \R \\ \eta (k) = 0 \, \text{ for  } k \in  \{ -1,0, N,N+1\} } } \frac{(\eta(k+1)- \eta(k))^2}{ 2 \, \Hscr_{[-1,N+1]}(\eta)} ,
\end{align*}
and to bound this quantity we write using the homogeneous boundary conditions
\begin{align*}
(\eta(k+1)- \eta(k))^2 &=\Big(  \sum_{j=0}^k  (\eta(j+1) - \eta(j)) - (\eta(j) - \eta(j-1))\Big)^2 \\
&\leq (k+1) \sum_{j=0}^k (\Delta \eta(j))^2 \leq N 2 \Hscr_{[-1,N+1]}(\eta),
\end{align*}
which establishes the desired bound on $\sigma(k)^2$. To derive the bound on $m(k)$ we make use of Lemma~\ref{le:Energytypical2}. First of all, by definition 
$m(k) = \phi^*(k+1) - \phi^*(k)$ where $\phi^*$ is the $\Hscr_{[-1,N+1]}$-minimiser subject to the boundary conditions $\br$ as well as the pinning condition $\Pscr$. 
If $(a, \alpha) \neq (0,0)$ and $(b,\beta) \neq (0,0)$ we invoke  Lemma~\ref{le:Energytypical2} for every fixed $M$ and for $N$ large enough (depending on $M$)  to get the existence of a profile $\tilde{\phi}^*$ satisfying the boundary conditions $\br$ as well as $\tilde{\phi}^*(p)=0$ for $(\delta_0 M)^{-\frac13} N \leq p \leq N - (\delta_0 M)^{-\frac13} N$ with $\Hscr_{[-1, N+1]}(\tilde \phi^*) \leq c \delta_0 MN$. In particular, this $\tilde{\phi}^*$ satisfies all of the $M$-typical pinning conditions simultaneously, which implies in turn that for each $M$-typical $\Pscr$ we have
 $\Hscr_{[-1,N+1]}(\phi^*) \leq \Hscr_{[-1,N+1]}(\tilde{\phi}^*)  \leq c \delta_0 MN$. Then we can write using the boundary conditions
\begin{align*}
|m(k)| =| \phi^*(k+1) - \phi^*(k)| &= \Big| \alpha N + \sum_{j=0}^k (\phi^*(j+1) - \phi^*(j)) - (\phi^*(j) - \phi^*(j-1) )\Big| \\
& \leq |\alpha| N + (k+1)^{\frac12} \Big(\sum_{j=0}^k (\Delta \phi^*)^2 \Big)^{\frac12} \leq   |\alpha| N  + 2c \delta_0 NM.
\end{align*}
At this point the required bound on $m(k)$ follows if we fix $\delta_0$ small enough. The argument if either $(a,\alpha) = (0,0)$ or $(b,\beta)=(0,0)$ is identical by noting that we can simply set $\tilde{\phi}^*$ near the corresponding boundary.

It remains to establish \eqref{e:Tight3} and to this end it suffices to derive an upper bound on 
\begin{align*}
\sum_{\Pscr \text{ not $M$-typcial}} \eps^{|\Pscr|}   Z_{[-1,N+1] \setminus \Pscr} (\br)
\end{align*}
as well as a lower bound on $\ZNe(\br)$. For the upper bound we fix $\Pscr$ and as before we denote by $\phi^*$ 
the unique  $\Hscr_{[-1,N+1]}$ minimiser subject to the boundary conditions $\br$ as well as the specific pinning condition $\Pscr$. 
As $\Pscr$ is not typical we can invoke Lemma~\ref{le:Energytypical1} to deduce that
$\Hscr_{[-1, N+1]}(\phi^*) \geq cN (\delta_0 M)^{\frac13}$. Then by Appendix~\ref{AppB} we have
 \begin{align*}
 Z_{[-1,N+1] \setminus \Pscr} (\br) = \ex^{-\Hscr_{[-1,N+1]}(\phi^*) }  Z_{[-1,N+1] \setminus \Pscr} (\zero) \leq \ex^{-cN (\delta_0 M)^{\frac13}} Z_{[-1,N+1] \setminus \Pscr} (\zero),
 \end{align*}
which permits to write  using \eqref{freeenergies}
\begin{align*}
\sum_{\Pscr \text{ not $M$-typcial}} \eps^{|\Pscr|}   Z_{[-1,N+1] \setminus \Pscr} (\br) \leq   \ex^{-cN (\delta_0 M)^{\frac13}} \sum_{\Pscr \subseteq \{ 1, \ldots, N-1\}} \eps^{|\Pscr|}   Z_{[-1,N+1] \setminus \Pscr} (\zero)
\leq  \ex^{-cN (\delta_0 M)^{\frac13}} \ex^{ N(\tau(\eps)  + \epsilon)}.
\end{align*}
On the other hand, for the lower bound we can use the coarse bound 
\begin{align*}
\ZNe(\br) =  \sum_{\Pscr \subseteq \{1, \ldots, N-1\}} \eps^{|\Pscr|}   Z_{[-1,N+1] \setminus \Pscr} (\br) \geq Z_{[-1,N+1]} (\br)=  \ex^{-NC}  Z_{[-1,N+1]}(\zero),
\end{align*}
where in the last step we have used that according to Proposition~\ref{minimiserconv} the mean energy $\frac{1}{N} \Hscr_{[-1,N+1]}$ is uniformly-in-$N$ bounded 
along the sequence of minimisers with boundary conditions $\br$ and without further pinning condition. This suffices to establish \eqref{e:Tight3}, 
 because according to Appendix~\ref{AppB} $Z_{[-1,N+1]}(\zero)$ decays at most polynomially.
\qed

\end{proofsect}

\begin{proofsect}{Proof of Lemma~\ref{le:Energytypical1}}
Assume first that both  $a \neq 0$ and $\alpha \neq 0$. Then for any $\phi$ satisfying the boundary condition $\phi(0) = N^2 a$ 
as well as $\phi(p) = 0$ for some $p \leq L$  there exists at least one $k \in \{0,1, \ldots, p-1 \}$
such that $|\phi(k+1) - \phi(k) | \geq  \frac{N^2 |a|}{p} \geq \frac{N^2 |a|}{L}$. We now recall that according to the boundary condition 
on $\phi(-1)$  we have $\phi(0) - \phi(-1) = N \alpha$. We now set $\delta_0 = \frac{a}{2\alpha}$ which implies that for 
$L \leq \delta_0 N$ we have $\frac{N^2 |a|}{L}  - N |\alpha| \geq  \frac12 \frac{N^2 |a|}{L}  $.
This then yields
\begin{align}
\notag
 \frac12 \frac{N^2 |a|}{L} &\leq  |\phi(k+1) - \phi(k) | -  |\phi(0) - \phi(-1)| \leq   |(\phi(k+1) - \phi(k) ) -  (\phi(0) - \phi(-1))|\\
 \notag
&=\Big| \sum_{j = 0}^k (\phi(j+1) - \phi(j) ) -  (\phi(j) - \phi(j-1))\Big| \leq (k+1)^{\frac12} \Big(\sum_{j=0}^k (\Delta \phi(j))^2 \Big)^{\frac12} \\
\label{e:Lemmm1}
& \leq L^{\frac12}  (2 \Hscr_{[-1, N+1]}(\phi))^{\frac12},
\end{align}
which can be rewritten as 
\begin{align*}
\Hscr_{[-1, N+1]}(\phi) \geq \frac{|a|^2}{8}\frac{N^4}{L^3},
\end{align*}
and which is stronger than the bound claimed in the proposition due to $L \leq N$.
If $\alpha =0$ and $a \neq 0$ the estimate \eqref{e:Lemmm1} holds in the same way without any restriction on $\delta_0$ and with left hand side replaced by  $  \frac{N^2 |a|}{L}$, 
i.e. the final lower bound on  $\Hscr_{[-1, N+1]}(\phi)$ is improved by an (irrelevant) factor $4$. Finally, let us assume $a=0$ and $\alpha \neq 0$, say $\alpha >0$. Then the  condition $\phi(p) = 0$ 
for some $p \leq L$ implies that there exists a $k \in \{0,1, \ldots, p-1 \}$ such that $\phi(k+1) - \phi(k) \leq 0$, so that \eqref{e:Lemmm1} holds with left hand side replaced by 
$N |\alpha|$ yielding the final estimate
\begin{align*}
\Hscr_{[-1, N+1]}(\phi) \geq \frac{|\alpha|^2 N^2}{L}.
\end{align*}

\qed
\end{proofsect}

\begin{proofsect}{Proof of Lemma~\ref{le:Energytypical2}}
%
%
%

We define the function 
$$
\tilde \phi^*(x)=\begin{cases} \phi^{*,\ell}(x) & \,\mbox{ for } x\in\{-1,0,1,\ldots,\floor{\delta N} , \floor{\delta N} +1\},\\
0 & \,\mbox{ for } x >L+1,\\
\phi^{*,r}(x) & \,\mbox{ for } x\in\{N - \floor{\delta N}-1,N- \floor{\delta N},\ldots,N , N +1\},
 \end{cases}
$$
where $ \phi^{*,\ell}$ is the minimiser for $\Hscr_{[-1,\floor{\delta N}+1]}(\phi)= \frac{1}{2}\sum_{k=0}^{\floor{\delta N}}(\Delta\phi_k)^2 $ (see Proposition~\ref{minimiserconv}) satisfying the boundary conditions 
\begin{equation*}
 \phi^{*,\ell}(-1)=N^2a - \alpha N,  \;\phi^{*,\ell}(0)=N^2a, \; \phi^{*,\ell}(\floor{\delta N})=0 , \; \mbox{ and }   \phi^{*,\ell}(\floor{\delta N})=0 ,
 \end{equation*}
 and similarly $\phi^{*,r}$ is the minimiser of $\Hscr_{[N -\floor{\delta N}-1, N+1]}(\phi)= \frac{1}{2}\sum_{k=N - \floor{\delta N} }^{N}(\Delta\phi_k)^2 $ satisfying
\begin{equation*}
 \phi^{*,r}(N - \floor{\delta N}-1)= 0;  \;\phi^{*,r}( N - \floor{\delta N})=0; \; \phi^{*,r}( N)=N^2 b , \; \mbox{ and }   \phi^{*,r}(N+1)=N^2 b + N \beta.
 \end{equation*}
Then $\Hscr_{[-1, N+1]}(\tilde \phi^*) = \Hscr_{[-1,\floor{\delta N}+1]}(\phi^{*, \ell}) + \Hscr_{[N -\floor{\delta N}-1, N+1]}(\phi^{*,r})$ and it remains to bound these two quantities. We only give the argument for 
 $\Hscr_{[-1,\floor{\delta N}+1]}(\phi^{*, \ell})$.


As in Proposition~\ref{minimiserconv} we  argue that
$$
\frac{1}{N}\Hscr_{[-1,\floor{\delta N}+1]}(\phi^{*, \ell}) \to \Escr_\delta(h^{*,(0,\delta)}_{(\ba,\zero)}) \quad\mbox{ as } N\to\infty,
$$
where $ h^{*,(0,\delta)}_{(\ba,\zero)} $ is the minimiser of $ \Escr_\delta(h)=\frac{1}{2}\int_0^\delta \ddot{h}^2(t)\,\d t $ with boundary conditions $ h(0)=a,\dot{h}(0)=\alpha,h(\delta)=0$ and $ \dot{h}(\delta)=0 $ (see Proposition~\ref{uniquemin}). Using Proposition~\ref{uniquemin} we compute
$$ 
\Escr_\delta(h^{*,(0,\delta)}_{(\ba,\zero)}) =\frac{1}{\delta^3}\big(6a^2+6a\alpha\delta +2\delta^2\alpha^2\big).
$$ Thus  for $N$ large enough
$$
\Hscr_{[-1,N+1]}(\phi)\le N\frac{2}{\delta^3} \big(6a^2+6a\alpha\delta +2\delta^2\alpha^2\big)
$$
as required.
\qed
\end{proofsect}

\begin{proofsect}{Proof of Lemma~\ref{lem:very-well-behaved}}
By sigma-additivity of the Lebesgue measure there exists a $\delta>0$ such  that
\begin{align*}
| \{ t \in [0,1] \colon |g(t) |<  2\delta \} | \leq  | \{ t \in [0,1] \colon g(t) = 0 \}| +\eps.
\end{align*}
For any $s \in  \{ t \in [0,1] \colon |g(t) |<  2\delta \} $ there exists a $\rho_s>0$ such that the ball $B(s, \rho_s) \cap [0,1]$ is still contained in this set. 
The collection of all these balls $ B(s, \rho_s) \cap [0,1]$  trivially covers $\{ t \in [0,1] \colon |g(t) |<  2\delta \} $ and therefore also the smaller and compact sub-level set
$\{ t \in [0,1] \colon |g(t) |\leq  \delta \}$. Thus there exists a finite collection $\{s_1, \ldots, s_{\tilde{M}} \}$ such that 
\begin{equation*}
\bigcup_{j=1}^{\tilde{M}} B(s_j, \rho_{s_j}) \supseteq \{ t \in [0,1] \colon |g(t) |\leq  \delta \} .
\end{equation*}
We then set 
$\mathcal{I} = \overline{\cup_{j=1}^{\tilde{M}} B(s_j, \rho_{s_j}) }  \cap [0,1]
$
and claim that this set has the desired properties. Indeed, the union of finitely many open intervals can always be written as the union 
of a (potentially smaller number of) \emph{disjoint} open intervals. The closure of such a set is the union of a finite (again, potentially smaller)
 number of disjoint closed intervals. The set $\mathcal{I}$ contains $\{ t \in [0,1] \colon |g(t) | \leq \delta \}$  by construction. Furthermore
\begin{equation*}
\bigcup_{j=1}^{\tilde{M}} B(s_j, \rho_{s_j}) \cap [0,1]  \subseteq \{ t \in [0,1] \colon |g(t) |<  2\delta \} 
\end{equation*} 
which implies that the measure of this set is bounded by $ | \{ s \in [0,1] \colon g(s) = 0 \}| +\epsilon$. Adding a finite number of boundary points 
does not change the Lebesgue measure, so that $\mathcal{I}$ satisfies the same bound.
\qed
\end{proofsect}
\end{proofsect}

\section*{Appendix} 
 \begin{appendices}
\section{Energy minimiser}\label{AppA}
 
 \noindent We outline the standard solution for the variational problem of minimising the \textit{energy functional}
 \begin{equation}\label{energy}
 \Escr(h)=\frac{1}{2}\int_0^1\,\ddot{h}^2(t)\,\d t\quad \mbox{ for } h\in H^2_{\br},
 \end{equation}
 where $ \br=(a,\alpha,b,\beta) $.
 
 \begin{proposition}\label{uniquemin}
 The variational problem, minimise $ \Escr$ in $ H^2_{\br} $, has a unique solution denoted by $ h_{\br}^*\in H^2_{\br}$ and given as
 $$
 h_{\br}^*(t)=a+\alpha t+k(\br)t^2+c(\br)t^3, \quad t\in [0,1],
 $$ with
 $$
 k(\br)=3(b-a)-2\alpha-\beta,\quad\mbox{ and } c(\br)=(\alpha+\beta)-2(b-a).
 $$
 Furthermore, $ \Escr(h_{\br}^*)=(2k(\br)^2+6k(\br)c(\br)+6c(\br)^2)$.
 \end{proposition}
 \begin{proofsect}{Proof}
For all $ h $ with $ h^{\ssup{4}}\equiv 0 $ we have  
 $$
 \langle \ddot{h},\ddot{g}-\ddot{h}\rangle_{L_2}=0 \quad\mbox{ for all } g\in H^2_{\zero}.
 $$ Then  we get
 $$
 \Escr(f)=\frac{1}{2}\langle \ddot{f},\ddot{f}\rangle_{L_2}\ge \frac{1}{2}\langle \ddot{h},\ddot{h}\rangle_{L_2}=\Escr(h).
 $$ We obtain the uniqueness by convexity and conclude with noting that $ (h^*_{\br})^{\ssup{4}}\equiv 0 $, see \cite{Mitrea} for an overview of bi-harmonic solutions.
 \qed
 \end{proofsect}

 \begin{lemma}\label{Gammaconv}
 For any $N \in \N$ let $h_N \colon \big\{-\frac{1}{N},0, \ldots, 1, \frac{N+1}{N} \big\} \to \R$ be given with boundary values $\br = ( a, \alpha, b, \beta)$ i.e.
 \begin{equation*}
 h_N(0) = a, \quad h_N(1) = b, \quad N\Big( h_{N}(0)- h_N\Big(-\frac{1}{N}\Big) \Big) = \alpha, \quad N \Big(h_N\Big(1+\frac{1}{N} \Big) - h_N(1) \Big) = \beta.
 \end{equation*}
We interpolate $h_N$ linearly between the grid-points. Furthermore set 
\begin{align*}
\Escr_N(h_N) = \frac12 \sum_{j=0}^N N^3 \Big[  h_N\Big(\frac{j+1}{N} \Big) +  h_N\Big(\frac{j-1}{N} \Big) -2 h_N\Big(\frac{j}{N} \Big)  \Big]^2.
\end{align*}
Then if $h_N$ converges uniformly over $[0,1]$ to a function $h$ we have 
\begin{align*}
\liminf_{N \to \infty}  \Escr_N(h_N) \geq \Escr(h)= \frac12 \int_0^1 \ddot{h}(t) ^2 \,\d t.
\end{align*}
\end{lemma}
 \begin{proofsect}{Proof}
We fix a subsequence  $(N_k)$ along which $\Escr_N(h_N)$ converges to $\liminf_{N \to \infty}  \Escr_N(h_N) $ which we can assume to be finite 
without loss of generality. Along this sequence $\Escr_N(h_N)$  is bounded. We drop the extra-index $k$ and assume from now on that
\begin{align}\label{newLemma1}
\sup_{N} \Escr_N(h_N) = \bar{C} < \infty.
\end{align}
We will first consider discrete derivatives of $h_N$. For $j=-1, \ldots, N$, we set
\begin{align}\label{newLemma2}
g_N \Big(\frac{j}{N} \Big) = N \Big[ h_N \Big(\frac{j+1}{N} \Big)   - h_N \Big(\frac{j}{N} \Big) \Big],
\end{align}
and as before we interpret $g_N$ as a function $[-\frac1N, 1] \to \R$ by linear interpolation between the 
grid-points. 
The functional $ \Escr_N(h_N) $ can be re-expressed in terms of $g_N$ as
\begin{align*}
\Escr_N(h_N) = \frac12 \sum_{j=0}^N N \Big[ g_N \Big(\frac{j}{N} \Big)   - g_N \Big(\frac{j-1}{N} \Big) \Big]^2 = \frac12 \int_{-\frac{1}{N}}^1 g_N'(s)^2 \,\d s.
\end{align*}
So, \eqref{newLemma1}  immediately implies the uniform H\"older bound
\begin{align}\label{newLemma3}
|g_N(t) - g_N(s) | \leq \int_{s}^t  |g_N'(r)| \,\d r \leq |t-s|^{\frac12}\Big( \int_{-{\frac1N}}^1 g_N'(r)^2 \,\d r\Big)^{\frac12} \leq |t-s|^{\frac12}  \sqrt{2 \bar{C}}
\end{align} 
for $-\frac{1}{N}\leq s<t \leq 1$. We introduce the (slightly) rescaled function $\tilde{g}_N \colon [0,1] \to \R$ defined 
as
\begin{align*}
\tilde{g}_N(t) = g_N \Big( \frac{N}{N+1} \Big( t+ \frac{1}{N} \Big) \Big)
\end{align*}
and observe that
\begin{align*}
 \int_0^1 \tilde{g}_N'(t)^2 dt = \frac{N}{N+1} \int_{-\frac{1}{N}}^1 g'_N(t)^2 \,\d t \leq 2 \bar{C}.
\end{align*}
Observing that $\tilde{g}_N(1)  = \beta$ we can conclude that there is a subsequence $N_k$ along which $\tilde{g}_{N_k}$ converges weakly in $H^1([0,1])$ to a function $g$ which 
satisfies
\begin{align*}
\int_0^1 g'(t)^2 \,\d t \leq \liminf_{N \to \infty}  \int \tilde{g}_{N_k}'(t)^2 \,\d t = \liminf_{N \to \infty} \frac{N}{N+1}  \int g_{N_k}'(t)^2 \,\d t = \liminf_{N \to \infty}\Escr(h_N).
\end{align*}
Thus, the desired statement follows as soon as we have established that for all $t \in [0,1]$,
$$
h(t) = a +\int_0^t g(s)\,\d s,
$$ because then we get  $\frac12\int_0^1 \ddot{h}(s)^2 \,\d s  = \frac12 \int_0^1 g'(s)^2\,\d s $.
To see this we rewrite the defining relation \eqref{newLemma2} of $g_N$ for any $N$ and any $t \in [\frac{j}{N},\frac{j+1}{N}]$, $j \geq 0$, as
\begin{align}\label{newLemma4}
h_N(t) = a + \sum_{k=0}^{j-1} \frac{1}{N} g_N(k) +\frac{ (Nt - j)}{N} g_N(j) = a + \int_0^t \tilde{g}_N(s) \,\d s + E_N,
\end{align}
where the error term $E_N$ satisfies
\begin{align*}
E_N \leq \Big| \int_0^t g_N(s) \,\d s - \int_0^s \tilde{g}_N(s)  \,\d s \Big| + \Big|  \int_0^{t}\Big( g_N(s) - g_N\Big( \frac{\floor{sN} }{N}  \Big)\Big)\,\d s \Big|.
\end{align*}
The definition of $\tilde{g}_N$ together with a uniform boundedness of $g_N$ in $L^1([0,1])$ imply that the first term converges to zero as $N\to \infty$ while the second term  can be seen to go to zero by the 
uniform H\"older bound \eqref{newLemma3}.  We can then conclude by going back to \eqref{newLemma4} and noting that on the one hand $h_N(t)$ converges to $h(t)$ by assumption and that on the other hand
the weak convergence of $\tilde{g}_N$ in $H^1([0,1])$ implies that $\int_0^t \tilde{g}_N(s)\,\d s$ converges to $\int_0^t g(s) \,\d s$.

\qed
 \end{proofsect}

 \begin{proposition}\label{minimiserconv}
 For any $ N\in\N, N>2 $, the variational problem, minimise $ \Hscr_{\L_N} $ in $\Omega_{\br}^N $ has a unique bi-harmonic solution $ \phi^*_{\br,N} \in \Omega_{\br}^N $ satisfying
 \begin{equation}\label{biharm}
 \begin{cases} \Delta^2\phi^*_{\br,N}(x)=0 & \mbox{ for } x\in\{1,\ldots, N-1\},\\
 \phi^*_{\br,N}(x)=\psi^{\ssup{N}}(x) & \mbox{ for } x\in\partial\L_N=\{-1,0,N,N+1\}.\end{cases}
 \end{equation}

Then the coefficients of the discrete polynomial $ h_{\br,N}^* $ given as $ h^*_{\br,N}(\xi):=\frac{1}{N^2}\phi^*_{\br,N}(\xi N) $ for $ \xi\in \{-\frac{1}{N},0,\frac{1}{N},\ldots,1,\frac{N+1}{N}\} $, converge as $ N\to\infty $ to the coefficients of the unique biharmonic function $ h_{\br}^*$.  Moreover,  
 $$
 \begin{aligned}
 \frac{1}{N}\Hscr_{\L_N}(\phi_{\br,N}^*)\longrightarrow \frac{1}{2}\Escr(h_{\br}^*) \quad\mbox { as }N\to\infty.
\end{aligned}
 $$
 \end{proposition}
 \begin{proofsect}{Proof}
 Similar to Proposition~\ref{uniquemin} one can show that the unique minimiser $ \phi_{\br,N}^* \in\Omega^N_{\br} $ is a polynomial of order three such that $ \phi_{\br,N}^*(N\xi)=N^2h_{\br,N}^*(\xi) $ with
 $$
h^*_{\br,N}(\xi)=a_N+\alpha_N(\br) t+k_N(\br)\xi^2+c_N(\br)\xi^3, \qquad \xi\in \{-\frac{1}{N},0,\frac{1}{N},\ldots,1,\frac{N+1}{N}\},
 $$
 with
 $$
 \begin{aligned}
 a_N&=a; \quad \alpha_N(\br)=\frac{2b-a(2+3N)+N(3b+\alpha(N+1)-\beta)}{(N+1)(N+2)};\\
 k_N(\br)&=N\frac{(-\alpha+\beta+N(3(b-a)-2\alpha-\beta))}{(N+1)(N+2)};\\
 c_N(\br)&=N^2\frac{(2(a-b)+\alpha+\beta)}{(N+1)(N+2)}.
 \end{aligned}
 $$
We observe that the coefficients of the polynomials $ h_{\br,N}^* $  converge to the ones of $ h_{\br}^*$, that is,  $ \alpha_N(\br)\to\alpha, k_N(\br)\to k(\br) $, and $ c_N(\br)\to c(\br) $ as $ N\to\infty $. The convergence of the minimal mean energy follows immediately with the established convergence of the polynomials. 
 
 \qed
 \end{proofsect}
 
\section{Partition function}\label{AppB}
We collect some known results about the partition function for the case with no pinning (see \cite{Borecki} and \cite{BS99}).  The partition function with zero boundary condition $ \br=\zero $ is
\begin{equation}
\begin{aligned}
Z_N(\zero)&=\int\,\ex^{-\Hscr_{[-1,N+1]}(\phi)}\,\prod_{k=1}^{N-1}\d\phi_k\prod_{k\in\{-1,0,N,N+1\}}\delta_0(\d\phi_k)
=\int_{\R^{N-1}}\,\ex^{-\frac{1}{2}\langle w, B_{N-1} w\rangle}\prod_{i=1}^{N-1}\d w_i\\&=\Big(\frac{(2\pi)^{N-1}}{\det(B_{N-1})}\Big)^{1/2}=\frac{(2\pi)^{\frac{N-1}{2}}}{\big(\frac{1}{12}(2+(N-1))^2(3+4(N-1)+(N-1)^2)\big)^{1/2}},
\end{aligned}
\end{equation}
where the matrix $ B_{N-1} $ reads as 
$$
\begin{pmatrix}
6 & -4 & 1 & 0 & \cdots & 0\\
-4 & 6 & -4 & 1 & 0 & \cdots\\
1 & -4 & 6 & -4 & 1 & \cdots \\
0 & 1 & -4 & 6 & -4 & \cdots\\
\cdots & \cdots & \cdots & \cdots & \cdots & \cdots\\
0 & \cdots &  \cdots  &  1 &  -4   & 6
\end{pmatrix}.
$$

We can easily obtain the following relation for the partition functions with given boundary $ \br $ (via $ \psi^{\ssup{N}} $) and zero boundary condition $ \br=\zero $ for models without pinning. 
 \begin{equation}\label{partition}
Z_N(\br) =\exp\big(-\Hscr_{\L_N}(\phi_{\br,N}^*)\big)Z_N(\zero).
\end{equation}

\end{appendices}

  \section*{Acknowledgments}
  S. Adams thanks J.D.~Deuschel for discussions on this problem as well T.~Funaki and H.~Sakagawa for discussions on Laplacian models and the warm hospitality during his visit in July 2015. 
  H. Weber was supported by  EPSRC First Grant EP/L018969/1 and a Royal Society  University Research Fellowship.

\end{document}